\documentclass[a4paper, 10pt]{article}

\usepackage{graphics,graphicx}
\usepackage{amssymb,amsmath}
\usepackage{fullpage}
\usepackage{enumitem}
\usepackage{xcolor}

\usepackage{boldline}

\usepackage[capitalize]{cleveref}

\usepackage{thmtools}

\newtheorem{theorem}{Theorem}

\newtheorem{lemma}{Lemma}
\newtheorem{conj}{Conjecture}
\newtheorem{problem}{Problem}

\newtheorem{claim}{}[theorem]

\def \no {\noindent}

\def \sm {\setminus}
\def \es {\emptyset}

\newcommand*\sq{\mathbin{\vcenter{\hbox{\rule{0.75ex}{1.0ex}}}}}

\newenvironment{proof}[1][]%
{\noindent {\setcounter{equation}{0}\it Proof.
	}{#1}{}}{\hfill$\Box$\vspace{2ex}}

\newenvironment{proof2}[1][]%
{\noindent {\setcounter{equation}{0}\it Proof.
	}{#1}{}}{$\sq$}

\def\longbox#1{\parbox{0.90\textwidth}{#1}}

\begin{document}
\title{Reconfiguration graph for vertex colorings of\\  $(P_2+P_3, C_4)$-free graphs}

\author{ M. Belavadi\thanks{Computer Science Unit, Indian Statistical
Institute, Chennai Centre, Chennai 600029, India.  Supported by ANRF-NPDF, Govt. of India (No.~PDF/2025/005208). ORCID: 0000-0002-3153-2339. } \and  T.~Karthick\thanks{Computer Science Unit, Indian Statistical
Institute, Chennai Centre, Chennai 600029, India. Email: karthick@isichennai.res.in. ORCID: 0000-0002-5950-9093. Corresponding author.} }
\maketitle

\begin{abstract}
For a graph $G$, let $\chi(G)$ denote the chromatic number of $G$. Given  a graph $G$, the \emph{reconfiguration graph for the $k$-colorings}  of $G$, denoted by  ${\cal R}_k(G)$, is the graph whose vertices are the $k$-colorings of $G$ and two $k$-colorings are joined by an edge if they differ on exactly one vertex of $G$. A graph $G$ is \emph{$k$-mixing} if ${\cal R}_k(G)$ is connected, and  is {\em recolorable} if it is $k$-mixing for all $k> \chi(G)$.    In this paper, we give a complete characterization of   $(P_2+P_3, C_4)$-free graphs that are recolorable. Moreover,  we show that if $G$ is a  recolorable $(P_2+P_3, C_4)$-free graph, then  for any $k >\chi(G)$,  the   diameter of ${\cal R}_k(G)$ is at most 2$n^{2}$. Furthermore, we show that if $G$ is a ($P_2+P_3, C_4$)-free graph  on $n$ vertices with degeneracy $\rho(G)$, then for all $k > \rho(G)+ 1$,   the diameter of ${\cal R}_k(G)$ is at most $O(n^2)$. This confirms a conjecture of Cereceda  for the class of ($P_2+P_3, C_4$)-free graphs.  These  results generalize some  known results available in the literature.   
\end{abstract}


\section{Introduction}

All our  graphs  in this paper are  simple, finite and undirected.  For missing notation and terminology, we refer to  \cite{BLS}.  Let $G$ be a graph   and let $k$ be a positive integer.  A \emph{$k$-coloring} of $G$ is a function $\phi$ from the vertex-set $V(G)$ to the set of colors  $\{1,2,\ldots, k\}$ such that  for any two adjacent vertices $u$ and $v$ in $V(G)$, we have $\phi(u)\neq \phi(v)$.   The \emph{chromatic number} of a graph $G$, denoted by $\chi(G)$, is the least  positive integer $k$ such that $G$ has a $k$-coloring.   A graph $G$ is \emph{$k$-chromatic} if $\chi(G)=k$, and it is \emph{$k$-colorable} if $\chi(G)\leq k$.

Given  a graph $G$, the \emph{reconfiguration graph for the $k$-colorings}  of $G$, denoted by  ${\cal R}_k(G)$, is the graph whose vertices are the $k$-colorings of $G$ and two colorings are joined by an edge if they differ on exactly one vertex of $G$.  There are several applications of reconfiguration graph for vertex colorings such as  frequency reassignment problem, maintaining a firewall in a changing network, characterizing Glauber dynamics, Markov chains in statistical physics, and planning motion  including 3D printing and robot movement, etc.  We refer to a  detailed survey   \cite{Heuvel-Survey} and the references therein for more information and applications.

\begin{figure}
\centering
 \includegraphics{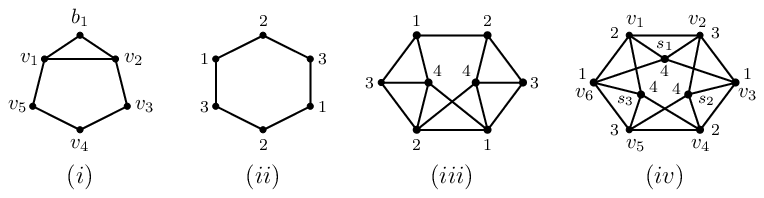}
\caption{Some special graphs:  $(i)$ A $5$-cap.  $(ii)$ The  graph $C_6$ and its frozen $3$-coloring. $(iii)$ The graph $F_1$ and its frozen $4$-coloring.
$(iv)$ The graph $F_2$ and its frozen $4$-coloring. }\label{cat}
\end{figure}

 A graph $G$ is \emph{$k$-mixing} if ${\cal R}_k(G)$ is connected and the \emph{$k$-recoloring diameter} of $G$ is the diameter of ${\cal R}_k(G)$. For instance,  ${\cal R}_{p}(K_{p})$ is disconnected, and thus $K_p$ is not $p$-mixing, where $K_p$ is the complete graph on $p$ vertices.   We say that a graph $G$ is {\em recolorable} if it is $k$-mixing for all $k> \chi(G)$.
A  $k$-coloring $\phi$ of a graph $G$ is called a  \emph{frozen $k$-coloring} (also called
\emph{fall $k$-coloring} \cite{Dunbar}), if for each $v\in V(G)$, all the $k$ colors appear in the closed neighborhood of $v$ (see \cref{cat} for examples). Note that  a  frozen $k$-coloring  of $G$   corresponds to an isolated vertex in  ${\cal R}_{k}(G)$.
Moreover, if a graph $G$ has a frozen $k$-coloring for some $k> \chi(G)$, then ${\cal R}_{k}(G)$ is disconnected, and hence $G$ is not $k$-mixing.

  Cereceda, van den Heuvel  and Johnson \cite{Cereceda2008} gave a polynomial-time algorithm to decide whether there
exists a path between two 3-colorings of a graph $G$ in ${\cal R}_3(G)$. However, Bonsma and Cereceda \cite{BC2009} proved that for all $k \geq 4$, it is \textsf{PSPACE}-complete to decide whether there exists a path between any two $k$-colorings of a graph $G$ in ${\cal R}_k(G)$.
For $k\in \{1,2\}$, it is easy to determine whether  ${\cal R}_k(G)$ is connected. But   it is \textsf{co-NP}-complete to decide   whether ${\cal R}_3(G)$ is connected     when $G$ is a bipartite graph \cite{CHJ2007}.

A class of graphs $\mathcal{G}$ is said to be \emph{hereditary} if it is closed under taking induced subgraphs. We say that a graph $G$ {\em contains} a graph $H$, if there is an induced subgraph of $G$ which is isomorphic to $H$.
A graph $G$ is   {\em $H$-free} if $G$ does not contain  $H$. For a set of graphs $\mathcal{H}$, a graph $G$ is  said to be {\em $\mathcal{H}$-free} if, for each $H$ in $\mathcal{H}$, $G$ is $H$-free. For any two vertex disjoint graphs $G$ and $H$, their {\em  union} (denoted by  $G + H$) is the  graph   with vertex-set $V(G)\cup V(H)$ and edge-set $E(G)\cup E(H)$, and their {\em join} (denoted by  $G\vee H$)  is the graph obtained from $G + H$ by adding an edge from each vertex of $G$ to every vertex of $H$.  The union of $t$   copies of $G$ is denoted by  $tG$.
For any positive integer $t$,  let $P_t$  and $C_t$ respectively denote the chordless path and  the chordless cycle
on $t$  vertices.

Structural properties of the reconfiguration graph of some hereditary classes of graphs have been studied in the past few years and gathered a wide attention recently;  see \cite{Manoj-thesis}.
Here, we are interested in  the following classification problem: ``{\em Given a hereditary class of graphs $\cal G$ and a graph $G\in \cal G$, classify whether $G$ is recolorable or not.}"  The class of perfect graphs and the class of graphs  defined by one or more forbidden induced subgraphs are among the most important and well-studied hereditary classes of graphs.  A graph is {\it perfect} if for each of its induced subgraph $H$, the chromatic number of $H$ equals the size of a largest complete subgraph of $H$.    The connectivity of the reconfiguration graph of  $k$-colorings for some classes of perfect graphs, such as bipartite graphs, chordal graphs, complement graphs of chordal graphs (or co-chordal graphs), weakly chordal graphs, $P_3$-free graphs, etc., were investigated very recently; see  \cite{ref3, Bonamy2018, Cereceda2008, FF-WeaklyChord}. Throughout this section,  unless stated otherwise, we use $\ell$ to denote a positive integer greater than the chromatic number of the given graph, and we use $n$ to denote the number of vertices in the given graph.

Cereceda, van den Heuvel  and Johnson \cite{Cereceda2008} showed that for all $\ell \ge 3$, there is a bipartite graph that is not $\ell$-mixing. Since the class of bipartite graphs is a subclass of the class of $K_3$-free graphs, it follows that an arbitrary   $K_3$-free graph   need not be recolorable. Also  for all $k\geq 4$, there is a $k$-colorable $2K_2$-free graph that admits a frozen $(k+1)$-coloring, and hence an arbitrary   $2K_2$-free graph need not be recolorable; see \cite{Manoj-thesis}. On the other hand, Bonamy and Bousquet \cite{Bonamy2018} showed that every $P_3$-free graph  is recolorable and the $\ell$-recolouring diameter is at most $2n$, and  that every $P_4$-free graph $G$ is recolorable and the $\ell$-recolouring diameter is at most $2n\chi(G)$. In \cite{ref3},  the first author with Cameron and Merkel  proved  that  if $H$ is an induced subgraph of   a $P_3+P_1$, then every $H$-free graph   is recolorable.  Moreover, if $H$ is any graph on at least five vertices, then there is an $H$-free graph which is not recolorable  \cite{ref3}. Thus we have the following problem.

 \begin{problem}[\cite{Manoj-thesis}]\label{q1}
Let $\cal{H}$ be a fixed set of graphs.   Given a $\cal{H}$-free graph $G$, classify whether $G$ is recolorable or not.
\end{problem}

As stated in the previous paragraph, Problem~\ref{q1} has been answered partially when ${\cal H}:=\{H\}$ and  $H$ has at most four vertices.
  However, when ${\cal H}:=\{H\}$ and  $H$ has at least five vertices, Problem~\ref{q1} is wide open. It is known that for all $\ell \ge 4$, there is a $P_5$-free graph that is not $\ell$-mixing, and recently   Lei et  al.  \cite{ChinesePaper} proved that   3-chromatic $P_5$-free graphs are recolorable. Below we state some results  when ${\cal H}$ has at least two graphs which are relevant to this paper.

A graph is a {\em split graph} if its vertex-set can be partitioned into a clique and an independent set. Note that a graph is a split graph if it is  both
chordal and co-chordal. From a result of \cite{Bonamy2018}, the class of split graphs is recolorable.
  An earlier known result of F\"oldes and Hammer \cite{FH-Split}  states that a graph is a split graph if and only if it is $(2K_2,C_4,C_5)$-free.
A graph is a {\em pseudo-split graph} if it is ($2K_2,C_4$)-free. From \cite{ref2}, it is known that    every ($2K_2,C_4$)-free graph is recolorable and the $\ell$-recolouring diameter is at most
$4n$. Belavadi et al. showed  that every ($P_5,C_4$)-free graph is recolorable  \cite{ref2}, and that every ($P_5,\overline{P_5},C_5$,~co-banner)-free graph is recolorable \cite{ref3}.  Feghali and Fiala \cite{FF-WeaklyChord} showed that every $3$-colorable $(P_5,\overline{P_5},C_5)$-free graph is recolorable.

    The {\em degeneracy} of a graph $G$, denoted by $\rho(G)$, is   the smallest positive integer $d$ such that  every subgraph of $G$ contains a vertex of degree at most $d$. Dyer et al. \cite{Dyer} (and independently  Cereceda et al. \cite{Cereceda2008}) showed  that if $G$ is any graph, then for each $k\geq \rho(G)+2$, ${\cal R}_k(G)$ is connected.  In 2007, Cereceda \cite{Cer-Thesis}  posed the following conjecture.

\begin{conj}[\cite{Cer-Thesis}]\label{Cer-conj}
If $G$ is a graph   on $n$ vertices, then for each $k\geq \rho(G)+2$,   the   diameter of  ${\cal R}_{k}(G)$ is  at most  $O(n^2)$.
\end{conj}

\cref{Cer-conj} is   open, and is known to be true  for trees \cite{BJLPP},   graphs $G$ with $\rho(G)=2$ and maximum degree at most $3$ \cite{FJP}, and for some special classes of graphs such as pseudo split graphs; see \cite{ref3}. In \cite{BH}, Bousquet and Heinrich proved that for any graph $G$ and for any $k\geq \rho(G)+2$, the diameter of ${\cal R}_{k}(G)$ is $O(n^{\rho(G)+1})$. See \cite{BH, Feghali2021} for  recent results related to   \cref{Cer-conj}.

In this paper, we focus  on \cref{q1} and \cref{Cer-conj} for a generalization  of split graphs and pseudo-split graphs, more specifically, for the class of ($P_2+P_3,C_4$)-free graphs. Indeed, we give a complete characterization of   $(P_2+P_3, C_4)$-free graphs that are recolorable. Moreover,  we show that if $G$ is a  recolorable $(P_2+P_3, C_4)$-free graph, then  for any $\ell >\chi(G)$,  the   diameter of ${\cal R}_\ell(G)$ is at most 2$n^{2}$.   Furthermore, we show that  \cref{Cer-conj}  holds for the class of ($P_2+P_3, C_4$)-free graphs.
In fact, the proofs of these results rely on a suitable structure theorem  for $(P_2+P_3, C_4)$-free graphs   that consists of some interesting special graph classes (see \cref{thm:str-c6-free,case-c6})  and  a few preliminary results on recoloring  which may be of independent interest.

\medskip
We finish this section with some notation and terminology used in this paper.

\medskip
\no{\bf Notation and terminology}: Given a graph $G$, we denote its  vertex-set by $V(G)$ and its edge-set by $E(G)$. For a graph $G$, we use $\overline{G}$ to denote the complement graph of $G$, and we use $n_{G}$ to denote the number of vertices in $G$; we simply use $n$ when the context is clear. For any $x, y\in V(G)$, we say that $x$ is a {\em neighbor} (resp. {\em  nonneighbor}) of   $y$ in $G$, if $x$  and $y$ are adjacent
(resp. nonadjacent) in $G$.

 Let $x$ be a vertex in $G$. The {\em neighborhood} of $x$, denoted by $N(x)$, is the  set of neighbors of $x$ in $G$, and  the {\em degree} of $x$  is the cardinality of $N(x)$.  The {\em closed neighborhood} of $x$, denoted by $N[x]$, is the set $\{x\}\cup N(x)$.

 For a subset $X$ of $V(G)$, we write $N(X)$ to denote the set $\{u\in V(G)\sm X \mid u  \mbox{ has a neighbor in }   X\}$. Two vertices $x$ and $y$ of a graph $G$ are said to be \emph{comparable} if they are nonadjacent, and $N(x)\subseteq N(y)$ or $N(y)\subseteq N(x)$. Given two disjoint subsets $X$ and $Y$ of $V(G)$,  we say that $X$ is {\em complete} to $Y$ if each vertex in $X$ is adjacent to every vertex in $Y$, and we say that   $X$ is {\em anticomplete} to $Y$ if no vertex in $X$ is adjacent to a vertex in $Y$.

Let $G$ be a graph. A {\em clique}  in $G$ is a set of mutually adjacent   vertices in $G$. The size of a maximum clique in $G$ is denoted by $\omega(G)$.
An {\em independent set}  in $G$ is a set of mutually nonadjacent vertices in $G$.


  A {\em blowup} $\cal B$ of a graph $G$ is the graph obtained from $G$ by replacing each vertex of $G$ with a  clique, and two cliques are complete to each other in $\cal B$ if  their corresponding vertices in $G$ are adjacent in $G$, and are anticomplete to each other in $\cal B$, otherwise.

  Given a positive integer $k$, we simply write $[k]$ to denote the set $\{1,2,\ldots, k\}$.

\medskip

From now on, the paper is organized as follows. 
We prove a structure theorem for the class of $(P_2+P_3, C_4)$-free graphs in \cref{sec:struc},  and  some useful preliminary results on recoloring in \cref{sec:pre-results}.  Finally, in \cref{sec:recol}, we   prove our main results (see \cref{recol-main,thm:degen}).

\section{Structure of ($P_2+P_3, C_4$)-free graphs}\label{sec:struc}

In this section, we give a structure of  ($P_2+P_3, C_4$)-free graphs which is suitable for proving our results on recoloring. Indeed, we   split our structure theorem into two cases based on whether the given  ($P_2+P_3, C_4$)-free graph contains a $C_6$ or not, and are given below in \cref{sec:C6-free,sec-H4}, respectively.

For reader's convenience, we give a schematic representation  for some of the new graph classes defined in this paper (see \cref{gcl2} and \cref{gcl})  and we use the following: A black dot represents a vertex and a shaded circle represents a clique. A solid line between any two shapes represents that the corresponding vertex-sets are complete to each other,  and the absence of a line between two shapes represents that the corresponding vertex-sets are anticomplete to each other.

 First we recall some structural properties of ($P_2+P_3, C_4$)-free graphs that contain a $C_5$ and use them in the latter sections (we  note that some of these properties are available in \cite{CK}).

\subsection{Properties of ($P_2+P_3, C_4$)-free graphs that contain a $C_5$} \label{genprop}

Let $G$ be a  ($P_2+P_3, C_4$)-free graph  that contains a $C_5$,  with vertex-set
$C:=\{v_1,v_2,v_3,v_4, v_5\}$ and edge-set $\{v_1v_2, v_2v_3, v_3v_4, v_4v_5, v_5v_1\}$.  For $i\in [5]$ and $i\mod 5$, we
let
 \begin{center}
\begin{tabular}{ccl}
$A_i$  &:= &  $\{u \in V(G)\setminus C\mid N(u) \cap C = \{v_i\}\}$,\\
$B_i$ &:= & $\{u \in V(G)\setminus C\mid N(u) \cap C =  \{v_i, v_{i+1}\}\}$,\\
$D_i$  &:= & $\{u \in V(G)\setminus C\mid N(u) \cap C =  \{v_{i-1}, v_i, v_{i+1}\}\}$,\\
$Z$ &:=& $\{u \in V(G)\setminus C\mid N(u) \cap C = C\}$, and\\
$T$ &:=& $\{u \in V(G)\setminus C\mid N(u) \cap C = \emptyset\}$.
\end{tabular}
\end{center}
Further we let $A:=\cup_{i=1}^5 A_i$, $B:= \cup_{i=1}^5
B_i$  and $D:= \cup_{i=1}^5 D_i$. Then since $G$ does not contain a $C_4$, we see that   $V(G)=C\cup  A\cup B\cup D\cup Z\cup T$. Moreover the following properties
hold where $i\in  [5]$ and $i\mod 5$: (We omit the proofs of  the properties which are available in \cite{CK},  and we give a proof for the new properties.)

\begin{enumerate}[label=  $(\mathbb{O}\arabic*)$, leftmargin=1cm] \itemsep=0pt

 \item\label{AiBi-card}
  $|A_i\cup B_i| \leq 1$ and  $|A_{i+1}\cup B_i| \leq 1$.

\item\label{AiAi+1}  $A_i$ is anticomplete to $A_{i-1}\cup A_{i+1}\cup D_{i-1}\cup D_{i+1}$.

\item \label{AD}   $A_i$ is complete to $D_i$. Moreover,  if
$A_i$ is nonempty, then $D_{i-2} \cup D_{i+2}$ is an empty set.

\item\label{DiDi+2-adj}
  $D_i$ is a clique,  and  $D_i$ is anticomplete to $D_{i+2}$.

\item\label{Z-adj}  $Z$ is a clique, and $Z$ is complete to   $C \cup D$.

\item\label{T-adj} $T$ is an independent set, and $T$ is anticomplete to  $A \cup B \cup D$.

\item\label{T-nemp}If $T$ is nonempty, then $G$ has a pair of comparable vertices.

{\em Proof}.~Let $t\in T$. Then from \ref{T-adj}, we have $N(t) \subseteq Z \subseteq N(v_1)$, and hence $t$ and $v_1$ are comparable vertices in $G$. This proves \ref{T-nemp}. $\sq$

\item\label{no-cap-BD}If $G$ does not contain a $5$-cap (see \cref{cat}:$(i)$), then $B$ is an empty set,  and  $D_i$ is complete to $D_{i+1}$  for each $i\in [5]$. Hence $C\cup D$ induces a blowup of a $C_5$.

{\em Proof}.~Suppose that $G$ does not contain a $5$-cap. If there is a vertex, say $b\in B$, then $C\cup \{b\}$ induces a $5$-cap which is a contradiction; so $B=\es$. Next if there are vertices, say $d_i\in D_i$ and $d_{i+1}\in D_{i+1}$, such that $d_id_{i+1}\notin E(G)$, then $(C\sm \{v_i\})\cup \{d_i,d_{i+1}\}$ induces a $5$-cap in $G$, a contradiction; so  $D_i$ is complete to $D_{i+1}$. Hence $C\cup D$ induces a blowup of a $C_5$, by \ref{DiDi+2-adj}. This proves  \ref{no-cap-BD}. $\sq$

\end{enumerate}

\subsection{Structure of ($P_2+P_3, C_4, C_6$)-free graphs}\label{sec:C6-free}
In this section, we study the structure of  ($P_2+P_3, C_4, C_6$)-free graphs. We begin with the following straightforward lemma.

\begin{lemma}\label{Lem:5-cap-free}
    Let $G$ be a connected ($P_2+P_3, C_4, C_6$, $5$-cap)-free graph that  contains a $C_5$. Then   $G$ has a pair of comparable vertices or for some $p\ge 0$, $G$ is the join of a  $K_p$ and a blowup of a $C_5$.
    \end{lemma}
    \begin{proof}
Suppose that $G$ contains a  $C_5$, say  with vertex-set  $C := \{v_1, v_2, v_3, v_4, v_5\}$ and edge-set $\{v_1v_2, v_2v_3, $ $ v_3v_4, v_4v_5,  v_5v_1\}$. Then with respect to $C$, we define the sets $A$, $B$, $D$, $Z$ and $T$ as in  \cref{genprop}, and we use the properties in  \cref{genprop}. Since $G$ is $5$-cap-free, we have $B=\es$ (by \ref{no-cap-BD}). Also   we may assume that $T = \es$ (otherwise, from  \ref{T-nemp}, $G$ has a pair of comparable vertices, and we are done). Now suppose that $A\neq \es$; let  $a_1 \in A_1$. Then $A_1=\{a_1\}$ (by \ref{AiBi-card}), and $D_3\cup D_4$ is an empty set (by \ref{AD}).       Further if there is a vertex in $A_{3} \cup A_{4}$, say $a'$, such that $a_1a'\in E(G)$, then
 $\{v_1,a_1,a',v_{3},v_{4},v_5\}$ or $\{v_1,v_2,v_3,v_4,a',a_1\}$ induces a $C_6$, a contradiction; so $\{a_1\}$ is anticomplete to $A_3\cup A_4$. Now  from \ref{AiBi-card}, \ref{AiAi+1}  and \ref{AD}, and by the above properties, we conclude that $N(a_1) \subseteq \{v_1\} \cup D_1\cup Z \subseteq N(v_2)$. Thus, $a_1$ and $v_2$ are comparable vertices in $G$, and we are done.   So suppose that $A=\es$.  Then from \ref{no-cap-BD}, we see that   $G[C\cup D]$ is a blowup of a $C_5$, and then from \ref{Z-adj}, we conclude that $G$ is the join of $G[C\cup D]$ (which is a blowup of a $C_5$) and $G[Z]\cong K_p$, where $p\ge 0$.  This proves \cref{Lem:5-cap-free}.
 \end{proof}

To proceed further, we define some new  graph classes.

\medskip

\begin{figure}

  \hspace{2cm}
 \includegraphics{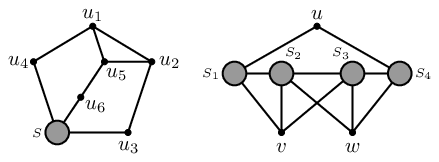}

\vspace{-3.1cm}

 \hspace{9.75cm}
 \includegraphics[height =3cm, width=3cm]{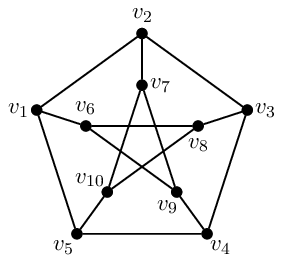}

 \caption{ Schematic representation of  ${\cal H}_1$ and ${\cal H}_2$, and the  Petersen graph  (left to right).}\label{gcl2}
  \end{figure}

\noindent{\bf Graph class ${\cal H}_1$}~(see \cref{gcl2}):~We say that a graph $G\in {\cal H}_1$  if its vertex-set can be partitioned into a set consisting of six vertices, say $\{u_1,u_2,\ldots, u_6\}$, and a nonempty clique, say $S$, such that the following hold:

$\bullet$~ $u_1u_2,u_2u_3, u_4u_1, u_1u_5,u_2u_5$ and $u_5u_6$ are edges in $G$.

$\bullet$~ $S$ is complete to $\{u_3,u_4,u_6\}$.

$\bullet$~ No other edges in $G$.

\medskip
\noindent{\bf Graph class ${\cal H}_2$}~(see \cref{gcl2}):~We say that a graph $G\in {\cal H}_2$  if its vertex-set can be partitioned into an independent set consisting of three vertices, say $\{u,v,w\}$, and four mutually disjoint cliques, say $S_1$, $S_2$, $S_3$ and $S_4$, such that the following hold:

$\bullet$~$S_1$, $S_2$ and $S_4$ are nonempty, and $S_3$ may be empty.

$\bullet$~$S_2$ is complete to $S_1\cup S_3$, and $S_3$ is complete to  $S_4$.

$\bullet$~$\{u\}$ is complete to $S_1\cup S_4$, $\{v\}$ is complete to $S_1\cup S_2\cup S_3$, and $\{w\}$ is complete to $S_2\cup S_3\cup S_4$.

$\bullet$~No other edges in $G$.

\begin{theorem}\label{thm:str-c6-free}
    Let $G$ be a connected ($P_2+P_3, C_4, C_6$)-free graph. Then at least one of the following hold:

          $\bullet$~$G$ is chordal.

      $\bullet$~$G$ has a pair of comparable vertices.

      $\bullet$~$G$ is the join of the graph $H$ and a $K_p$, where $H$ is a blowup of a $C_5$ or $H\in \mathcal{H}_1 \cup \mathcal{H}_2$, and $p\ge 0$.

   \end{theorem}
\begin{proof}
Let $G$ be a connected ($P_2+P_3, C_4, C_6$)-free graph. We may assume that $G$ is not chordal and that $G$ does not have a pair of  comparable vertices. Now since $G$ is ($P_2+P_3, C_4, C_6$)-free and since for each $k\geq 7$, the graph $C_k$  contains a $P_2+P_3$, we may assume that $G$ contains a $C_5$.  If $G$ is $5$-cap-free, then the result follows from \cref{Lem:5-cap-free}. So suppose that $G$ contains a $5$-cap. We may assume such a $5$-cap in $G$  with vertices and edges as shown in \cref{cat}, and we let $C := \{v_1, v_2, v_3, v_4,$ $ v_5\}$. Then with respect to $C$, we define the sets $A$, $B$, $D$, $Z$ and $T$ as in  \cref{genprop}, and we use the properties in  \cref{genprop}. Clearly $b_1 \in B_1$, and so $B_1=\{b_1\}$, by \ref{AiBi-card}. By \ref{T-nemp}, we may assume that $T=\es$.   Also the following hold.

      \begin{claim}\label{AB-emp}   $A\sm A_4 $ and $B\sm \{b_1\}$ are empty sets. So, $V(G)=C\cup A_4\cup \{b_1\}\cup D\cup Z$.
      \end{claim}
      \begin{proof2}
      Since $B_1=\{b_1\}$,  by using \ref{AiBi-card},  it is enough to show that $A_3\cup A_5\cup (B\sm \{b_1\}) =\es$. By   symmetry, we show that
$A_3\cup  B_2\cup B_3=\es$. Suppose to the contrary that there is a vertex, say $u\in A_3\cup  B_2\cup B_3$. Now if $ub_1\notin E(G)$, then $\{u,v_3, b_1,v_1,v_5\}$  induces a $P_2+P_3$; so  $ub_1\in E(G)$. Thus if $u\in A_3\cup B_3$, then $\{u,b_1,v_2,v_3\}$ induces a $C_4$; so $u\in B_2$, and then $\{v_1,b_1,u,v_3,v_4,v_5\}$ induces a $C_6$ which is a contradiction. This proves \ref{AB-emp}.
\end{proof2}

      \begin{claim}\label{BD-adj}
        $\{b_1\}$ is complete to $D\sm D_4$,  and  $\{b_1\}$ is anticomplete to $D_4$.
      \end{claim}
 \begin{proof2}
 By   symmetry, we show that $\{b_1\}$ is complete to $D_1\cup D_3$, and $\{b_1\}$ is anticomplete to $D_4$.  If there is a vertex, say
 $d\in D_1\cup D_3$, such that $b_1d\notin E(G)$, then $\{v_3,v_4,b_1,v_1,d\}$ or $\{d,v_3,b_1,v_1,v_5\}$ induces a $P_2+P_3$ which is a contradiction; so $\{b_1\}$ is complete to $D_1\cup D_3$. Next if there is a vertex in $D_4$, say
 $d_4$, such that $b_1d_4\in E(G)$, then $\{b_1,v_1,v_5,d_4\}$ induces a $C_4$ which is a contradiction; so $\{b_1\}$ is anticomplete to $D_4$. This proves \ref{BD-adj}.
 \end{proof2}

   \begin{claim}\label{D-adj} For each $i\in
[5]$ and $i$ mod $5$, $D_i$ is complete to $D_{i+1}$.
   \end{claim}
   \begin{proof2}
   By  symmetry, we prove the assertion for $i\in \{1,2,3\}$. Suppose to the contrary that there are nonadjacent vertices, say $d_i\in D_i$ and $d_{i+1}\in D_{i+1}$, where  $i\in \{1,2,3\}$,  and we use \ref{BD-adj}. Now if $i=1$, then $\{d_1,b_1,d_2, v_3,v_4,v_5\}$ induces a $C_6$; if $i=2$, then
$\{b_1,d_2,v_3,d_3\}$ induces a $C_4$; and if $i=3$, then $\{v_1,b_1,d_3,$ $v_3,d_4,v_5\}$ induces a $C_6$ which are contradictions. This proves \ref{D-adj}.
   \end{proof2}

    \begin{claim}\label{D3D5-emp}  One of $D_{3}$ and $D_{5}$ is empty.
   \end{claim}
   \begin{proof2}
   If there are vertices, say $d_3\in D_3$ and $d_5\in D_5$, then $d_3d_5\notin E(G)$ (by \ref{DiDi+2-adj}), and then from \ref{BD-adj}, we see that $\{b_1,d_3,v_4,d_5\}$ induces a $C_4$ which is a contradiction. So \ref{D3D5-emp} holds.
   \end{proof2}

   \begin{claim}\label{A4B1} $A_4$ is complete to $\{b_1\}$.
\end{claim}
 \begin{proof2}
      If there is a vertex, say $a_4\in A_4$, such that $a_4b_1\notin E(G)$,  then $\{b_1,v_2, a_4,v_4,v_5\}$ induces a $P_2+P_3$ which is a contradiction. So \ref{A4B1} holds.
      \end{proof2}

    \begin{claim}\label{Z-A4B1} $Z$ is complete to $A_4\cup \{b_1\}$.
   \end{claim}
   \begin{proof2}
   First suppose that $A_4=\es$.  Suppose to the contrary that there is a vertex in $Z$, say $z$, such that $zb_1\notin E(G)$.
    Then by \ref{AB-emp}, \ref{BD-adj} and \ref{Z-adj}, we see that $N(b_1) = \{v_1,v_2\}\cup (D\sm D_4)\cup (Z\sm \{z\})\subseteq N(z)$, and hence $z$ and $b_1$ are comparable vertices in $G$ which is a contradiction to our assumption; so $Z$ is complete to $\{b_1\}$, and we  are done.
    So suppose that $A_4\neq \es$, and we let $A_4:=\{a_4\}$, by \ref{AiBi-card}. Then $a_4b_1\in E(G)$, by \ref{A4B1}. Let $z^*\in Z$ be  arbitrary.   Now if $z^*$ is anticomplete to $\{b_1,a_4\}$, then  $\{a_4,b_1,v_3,z^*,v_5\}$ induces a $P_2+P_3$ which is a contradiction; so  $z^*$ is adjacent to   one of $b_1$ and $a_4$. Then since exactly one of $\{z^*,a_4,b_1,v_2\}$ and $\{z^*,b_1,a_4,v_4\}$ does not induce a  $C_4$, we conclude that $z^*$ is complete to $\{a_4,b_1\}$, and we are done.  This proves \ref{Z-A4B1}.
   \end{proof2}

    \begin{claim}\label{A4D} If $A_4$ is nonempty, then $D\sm D_4$ is an empty set.
   \end{claim}
   \begin{proof2}Let $A_4:=\{a_4\}$, by \ref{AiBi-card}. By \ref{AD} and by  symmetry, it is enough to show that $D_3=\es$.  From  \ref{A4B1}, we have $a_4b_1\in E(G)$. Now if there is a vertex, say $d_3\in D_3$, then we have $a_4d_3\notin E(G)$ (by \ref{AiAi+1}), and $b_1d_3\in E(G)$ (by \ref{BD-adj}), and then $\{b_1,d_3,v_4,a_4\}$ induces a $C_4$ which is a contradiction; so $D_3=\es$. This proves \ref{A4D}.
   \end{proof2}

\medskip

      \begin{claim}\label{H-H1H2} Let $H:=G[C\cup A_4 \cup \{b_1\}\cup D]$. Then    $H$ is in $\mathcal{H}_1 \cup \mathcal{H}_2$.
   \end{claim}
   \begin{proof2}
     First suppose that $A_4\neq \es$; we let $A_4:=\{a_4\}$, by \ref{AiBi-card}. Then from \ref{A4D}, we have $D\sm D_4=\es$. So from \ref{AD}, \ref{DiDi+2-adj}, \cref{AB-emp},  \cref{BD-adj} and \cref{A4B1}, we conclude that $H\in {\cal H}_1$, where we take $u_1:=v_1, u_2:=v_2, u_3:=v_3, u_4:=v_5, u_5:=b_1, u_6:=a_4$ and $S:=\{v_4\}\cup D_4$, and we are done.  So we may  assume that $A_4 = \es$. From \ref{D3D5-emp}, we know that one of $D_3$ and $D_5$ is empty. By symmetry, we may assume that $D_5 = \es$. Then, from  \ref{DiDi+2-adj}, \ref{BD-adj} and \ref{D-adj}, we  see that $H\in \mathcal{H}_2$, where we take $u:=v_5$, $v:=b_1$, $w:=v_3$, $S_1:= D_1\cup \{v_1\}$, $S_2:= D_2 \cup \{v_2\}$, $S_3:= D_3$ and $S_4:= D_4 \cup \{v_4\}$, and we are done.
This proves \ref{H-H1H2}.
   \end{proof2}

\medskip
 From \ref{Z-adj} and \ref{Z-A4B1}, $Z$ is a clique which is complete to $C\cup A_4\cup \{b_1\}\cup D$,  and from \ref{AB-emp}  and \ref{H-H1H2}, we conclude that $G$ is the join of $H:=G[C\cup A_4\cup \{b_1\}\cup D]$ and $G[Z]\cong K_p$, where $H\in \mathcal{H}_1 \cup \mathcal{H}_2$  and $p \ge 0$.  This completes the proof.
 \end{proof}

 \begin{figure}
\centering
 \includegraphics{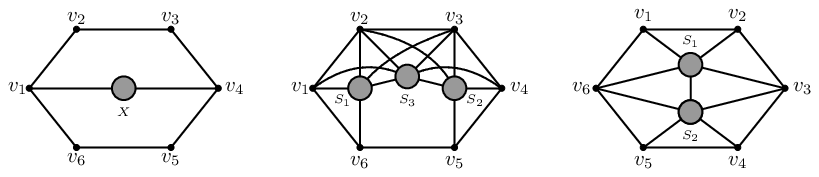}
\caption{ Schematic representation of ${\cal H}_3$, ${\cal H}_4$ and ${\cal H}_5$ (left to right).}\label{gcl}
\end{figure}

\subsection{Structure of ($P_2+P_3, C_4$)-free graphs that contain a $C_6$} \label{sec-H4}

In this section, we study the structure of  ($P_2+P_3, C_4$)-free graphs that contain a $C_6$. To proceed further,  we define a few more new  graph classes.

\medskip
\noindent{\bf Graph class ${\cal H}_3$}~(see \cref{gcl}):~We say that a graph $G\in {\cal H}_3$  if its vertex-set can be partitioned into a set consisting of six vertices, say $\{v_1,v_2,\ldots, v_6\}$, and a nonempty clique, say $X$, such that the following hold:

$\bullet$~$G[\{v_1,v_2,\ldots, v_6\}]\cong C_6$ with edge-set $\{v_1v_2, v_2v_3, v_3v_4, v_4v_5, v_5v_6,v_6v_1\}$.

$\bullet$~$X$ is complete to $\{v_1,v_4\}$.

$\bullet$~No other edges in $G$.

\smallskip
 \noindent The {\em theta graph}  is a graph in ${\cal H}_3$ with $|X|=1$.

\medskip
\noindent{\bf Graph class ${\cal H}_4$}~(see \cref{gcl}):~We say that a graph $G\in {\cal H}_4$  if its vertex-set can be partitioned into a set consisting of six vertices, say $\{v_1,v_2,\ldots, v_6\}$, and three mutually disjoint cliques, say $S_1$, $S_2$ and $S_3$,  such that the following hold:

$\bullet$~$G[\{v_1,v_2,\ldots, v_6\}]\cong C_6$ with edge-set $\{v_1v_2, v_2v_3, v_3v_4, v_4v_5, v_5v_6,v_6v_1\}$.

$\bullet$~$S_1$ is nonempty, and $S_2$ and $S_3$ may be empty.

 $\bullet$~$S_1$ is complete to $\{v_1,v_2,v_3,v_6\}$, $S_2$ is complete to $\{v_2,v_3,v_4,v_5\}$, and  $S_3$ is complete to $\{v_1,v_2,v_3,v_4\}$.

$\bullet$~$S_3$ is complete to  $S_1\cup S_2$.

$\bullet$~No other edges in $G$.

\smallskip
 \noindent We say that a graph $G\in {\cal H}_4^*$ if  $G\in {\cal H}_4$ with $|S_1|= |S_2|= p$, where $p\geq 1$ and $S_3$ is an empty set. Note that the graph $F_1$ (see \cref{cat}) belongs to ${\cal H}_4^*$ where $p=1$.

\medskip

\noindent{\bf Graph class ${\cal H}_5$}~(see \cref{gcl}):~We say that a graph $G\in {\cal H}_5$   if its vertex-set can be partitioned into a set consisting of six vertices, say $\{v_1,v_2,\ldots, v_6\}$, and two mutually disjoint nonempty cliques, say $S_1$ and $S_2$,  such that:

$\bullet$~$G[\{v_1,v_2,\ldots, v_6\}]\cong C_6$ with edge-set $\{v_1v_2, v_2v_3, v_3v_4, v_4v_5, v_5v_6,v_6v_1\}$.

 $\bullet$~$S_1$ is complete to $S_2\cup \{v_1,v_2,v_3,v_6\}$, and $S_2$ is complete to $\{v_3,v_4,v_5,v_6\}$.

$\bullet$~No other edges in $G$.

\medskip
We will use $\cal L$ to denote the set of   graphs $H^*$, where $H^*$  is an induced subgraph of the Petersen graph, and $H^*$   contain a theta graph.

\begin{theorem}\label{case-c6}
    Let $G$ be a connected ($P_2+P_3, C_4$)-free graph that contains a $C_6$. Then  $G$ has a pair of comparable vertices or
 $G$ is the join of ($H+\overline{K_q}$) and a $K_p$, where $H\in \{C_6, F_2\} \cup {\cal L} \cup \mathcal{H}_3\cup \mathcal{H}_4 \cup \mathcal{H}_5$,  $p\ge 0$ and $q\ge 0$.

\end{theorem}

\begin{proof} Let $G$ be a connected ($P_2+P_3, C_4$)-free graph  that contains a $C_6$, say with vertex-set
$C:=\{v_1,v_2,v_3,v_4, v_5,v_6\}$ and edge-set $\{v_1v_2, v_2v_3, v_3v_4, v_4v_5, v_5v_6,v_6v_1\}$.  For   $i\in [6]$ and $i$ mod $6$, we let:
 \begin{center}
\begin{tabular}{ccl}
$X_i$  &:= &  $\{u \in V(G)\setminus C\mid N(u) \cap C = \{v_i,v_{i+3}\}\}$,\\
$Y_i$ &:= & $\{u \in V(G)\setminus C\mid N(u) \cap C =  C\sm \{v_{i-1},v_{i-2}\}\}$,\\
 $Z$ &:=& $\{u \in V(G)\setminus C\mid N(u) \cap C = C\}$,  and\\
$T$ &:=& $\{u \in V(G)\setminus C\mid N(u) \cap C = \emptyset\}$.
\end{tabular}
\end{center}
Note that $X_i=X_{i+3}$ for each  $i\in [6]$ and $i$ mod $6$. Further we let $X:=\cup_{i=1}^3 X_i$ and $Y:= \cup_{i=1}^6
Y_i$. Throughout this theorem, the arithmetic operations on the indices are understood to be modulo 6.  We may assume that $G$ does not have a pair of  comparable vertices. First we show that:

\begin{claim}
 \label{C6-nei}   $V(G)=C\cup  X\cup Y\cup  Z\cup T$.
\end{claim}
\begin{proof2}
 Let $u\in V(G)\sm C$ be arbitrary. We will show that if $u\notin Z\cup T$, then $u\in X\cup Y$.  Since $u\notin Z\cup T$,  there is an index $i\in [6]$ such that $uv_i\in E(G)$ and $uv_{i-1}\notin E(G)$. By symmetry, we may assume that $i=1$. Then since $\{v_1,u,v_5,v_6\}$ does not induce a $C_4$, we have $uv_5\notin E(G)$. Now if $uv_2\in E(G)$, then since $\{u,v_2,v_4,v_5,v_6\}$ does not induce a $P_2+P_3$, we have $uv_4\in E(G)$, and then since $\{u,v_2,v_3,v_4\}$ does not induce a $C_4$, we must have $uv_3\in E(G)$  and hence $u\in Y_1$; so we may assume that $uv_2\notin E(G)$. Then since $\{u,v_1,v_2,v_3\}$ does not induce a $C_4$, we have $uv_3\notin E(G)$, and then since $\{u,v_1,v_3,v_4,v_5\}$ does not induce a $P_2+P_3$, we must have $uv_4\in E(G)$ which implies that $u\in X_1$. This proves \ref{C6-nei}.
\end{proof2}

\medskip
Moreover, the following hold:

\begin{claim}\label{C6-XYZT} For each $i\in [6]$,  $X_i \cup Y_i\cup Z$ is a clique.
\end{claim}
\begin{proof2}
 If there are nonadjacent vertices in $X_i \cup Y_i\cup Z$, say $p$ and $q$, then $\{v_i,p,v_{i+3},q\}$ induces a $C_4$ which is a contradiction; so \ref{C6-XYZT} holds.
\end{proof2}

\begin{claim}\label{C6-Xi}For $i\in [3]$, $X_i$'s are mutually anticomplete to each other.
\end{claim}
\begin{proof2}
 We show that $X_i$ is anticomplete to $X_{i+1}$, and then the proof follows from symmetry. Now if there are vertices, say $x_i\in X_i$ and $x_{i+1}\in X_{i+1}$, such that $x_ix_{i+1}\in E(G)$, then  $\{x_i, v_i, v_{i+1},x_{i+1}\}$ induces a $C_4$ which is a contradiction. So \ref{C6-Xi} holds.
\end{proof2}

\begin{claim}\label{C6-Xi2} For each $i\in [3]$, if $|X_i|\geq 2$, then $X_{i+1}\cup X_{i+2}$ is an empty set.
\end{claim}
\begin{proof2}
Let $x_i, x_i'\in X_i$. Then from \ref{C6-XYZT}, we have $x_ix_i'\in E(G)$. Now if  there is a vertex in $X_{i+1}$,  say $x_{i+1}$, then from \ref{C6-Xi}, we have $x_ix_{i+1}, x_i'x_{i+1}\notin E(G)$, and then $\{x_i, x_i', x_{i+1}, v_{i+1}, v_{i+2}\}$ induces a $P_2+P_3$ which is a contradiction;  so $X_{i+1}=\es$. Likewise,    $X_{i+2}=\es$. This proves \ref{C6-Xi2}.
\end{proof2}

\begin{claim}
\label{C6-Y}For each $i\in [6]$, $Y_i$ is complete to $Y_{i-1}\cup Y_{i+1}\cup Y_{i+3}$, and is anticomplete to $Y_{i-2}\cup Y_{i+2}$.
\end{claim}
\begin{proof2}
 If there are vertices, say $y\in Y_i$ and  $y'\in Y_{i-1}\cup Y_{i+1}\cup Y_{i+3}$, such that $yy'\notin E(G)$, then one of $\{y, v_i, y', v_{i+2}\}$,
 $\{y, v_{i+1}, y', v_{i+3}\}$ and $\{y, v_i, y', v_{i+3}\}$ induces a $C_4$ which is a contradiction; so  $Y_i$ is complete to $Y_{i-1}\cup Y_{i+1}\cup Y_{i+3}$. Also if there are vertices, say $y\in Y_i$ and $y^*\in Y_{i-2}\cup Y_{i+2}$, such that $yy^*\in E(G)$, then  $\{y,y^*, v_{i-2}, v_{i+3}\}$ or $\{y, y^*, v_{i-1}, v_{i}\}$
  induces a $C_4$ which is a contradiction; so  $Y_i$ is anticomplete to $Y_{i-2}\cup Y_{i+2}$. This proves \ref{C6-Y}.
\end{proof2}

\begin{claim}
\label{C6-Yi+3}For each $i\in [6]$, if  $Y_i$ and $Y_{i+1}$ are nonempty,  then $Y_{i+3} \cup Y_{i-2}$ is an empty set, and hence one of $Y_{i+2}$ and $Y_{i-1}$ is  empty.
\end{claim}
\begin{proof2}
 Let $y_i\in Y_i$ and $y_{i+1}\in Y_{i+1}$. Suppose to the contrary   that there is a vertex in $Y_{i+3}$, say $y_{i+3}$. Then from \ref{C6-Y}, we have $y_iy_{i+1}, y_iy_{i+3}\in E(G)$ and $y_{i+1}y_{i+3}\notin E(G)$. But then $\{y_i,y_{i+1},v_{i-2},y_{i+3}\}$ induces a $C_4$ which is a contradiction; so $Y_{i+3}=\es$. Likewise, $Y_{i-2}=\es$.
\end{proof2}

\begin{claim}
\label{C6-Yigeq2}For each $i\in [6]$, if $|Y_i|\geq 2$, then  one of $Y_{i-2}$ and $Y_{i+2}$ is  empty.
\end{claim}
\begin{proof2}
Let $y_{i},y_{i}'\in Y_i$.  Suppose to the contrary   that there are vertices, say $y_{i-2}\in Y_{i-2}$ and $y_{i+2}\in Y_{i+2}$. Then from \ref{C6-XYZT}, we have $y_{i}y_{i}'\in E(G)$ and from \ref{C6-Y}, we have $\{y_{i}, y_{i}'\}$ is anticomplete to $\{y_{i-2}, y_{i+2}\}$ and $y_{i+2}y_{i-2}\notin E(G)$. But then $\{y_{i},y_{i}',y_{i+2},v_{i-1},y_{i-2}\}$ induces a $P_2+P_3$ which is a contradiction.
\end{proof2}

\begin{claim}
\label{C6-XorY} One of $X$ and $Y$ is an empty set.
\end{claim}
\begin{proof2}
Suppose that $Y\neq \es$. We may assume that $y\in Y_1$. We show that $X=\es$. By symmetry, it is enough to show that  $X_1\cup X_2=\es$. Suppose to the contrary that there is a vertex in $X_1\cup X_2$, say $x$. If $x\in X_1$, then from \ref{C6-XYZT}, we have $yx\in E(G)$ and then $\{v_5,v_6,x,y,v_3\}$ induces a $P_2+P_3$ which is a contradiction; so we may assume that $x\in X_2$. Then we see that $\{y,x,v_5,v_4\}$ induces a $C_4$ or $\{y,v_3,x,v_5,v_6\}$ induces a $P_2+P_3$ which is a contradiction. This proves \ref{C6-XorY}.
\end{proof2}

\begin{claim}\label{C6-TY}
\label{C6-T} $T$ is an independent set, and $T$ is anticomplete to $Y$.
\end{claim}
\begin{proof2}
Suppose to the contrary that there are adjacent vertices, say $t\in T$ and $p\in T\cup Y$. By symmetry, we may assume that $p\in T\cup Y_1$.  Then $\{t,p, v_1,v_2,v_3\}$ or
  $\{v_5,v_6,t,p,v_3\}$ induces a $P_2+P_3$ which is a contradiction. So \ref{C6-TY} holds.
\end{proof2}

\begin{claim}\label{C6-T-X} If a vertex in $T$ has a neighbor in $X$, then it is complete to $X$.
\end{claim}
\begin{proof2}
Let $t$ be a vertex in $T$ that has a neighbor in $X$, say $x^*$. By symmetry, we may assume that $x^*\in X_1$. Suppose to the contrary that there is a vertex in $X$, say $x$, such that $tx\notin E(G)$. Then again by  symmetry, we may assume that $x\in X_1\cup X_2$. Now if $x\in X_1$, then from \ref{C6-XYZT}, we have $x^*x\in E(G)$, and then $\{v_2,v_3,t, x^*,x\}$ induces a $P_2+P_3$; so $x\in X_2$. Then from \ref{C6-Xi}, we have $x^*x\notin E(G)$, and then $\{t,x^*,x,v_5,v_6\}$ induces a $P_2+P_3$ which is a contradiction. So $\{t\}$ is complete to $X$. This proves \ref{C6-T-X}.
\end{proof2}

\begin{claim}\label{C6-T-Z}  $T$ is complete to $Z$.
\end{claim}
\begin{proof2}
If there are nonadjacent vertices, say $t\in T$ and $z\in Z$,  then from \ref{C6-nei}, \ref{C6-XYZT} and \ref{C6-TY}, we see that $N(t)\subseteq X \cup (Z\setminus\{z\}) \subseteq N(z)$ and hence $t$ and $z$ are comparable vertices in $G$ which is a contradiction to our assumption; so  $T$ is complete to $Z$. This proves \ref{C6-T-Z}.
\end{proof2}

\begin{claim}\label{C6-TZX}
    $T$ is complete to $X$.
\end{claim}
\begin{proof2}
    If there is a vertex in $T$, say $t$, which is anticomplete to $X$, then from \ref{C6-TY} and \ref{C6-T-Z}, we have $N(t)= Z \subseteq N(v_1)$, and hence $t$ and $v_1$ are comparable vertices in $G$ which is a contradiction to our assumption; so  each vertex in $T$ has a neighbor in $X$. Then it follows from \ref{C6-T-X} that $T$ is complete to $X$. This proves \ref{C6-TZX}.
\end{proof2}

\begin{claim}\label{C6-T-Xi}
    If there is an  $i\in [3]$ such that $X\sm X_i$ is an empty set, then   $T$ is an empty set.
\end{claim}
\begin{proof2}
    By symmetry, we may assume that $i=1$. Now if there is a vertex in $T$, say $t$, then from  \ref{C6-TY}, \ref{C6-T-Z} and \ref{C6-TZX}, we see that $N(t) = Z \cup X_1 \subseteq N(v_1)$ and thus $t$ and $v_1$ are comparable vertices in $G$ which is a contradiction to our assumption; so   $T=\es$. This proves \ref{C6-T-Xi}.
\end{proof2}

\begin{claim}\label{C6-Yi-Yi+1}
If there is an  $i\in [6]$ such that $Y_i$ and $Y_{i+1}$ are nonempty, then $G[C\cup Y] \in {\cal H}_4$.
\end{claim}
\begin{proof2}
 By symmetry, we may assume that $Y_6$ and $Y_1$ are nonempty. Then from \ref{C6-Yi+3} and again by symmetry, we may assume that $Y_3\cup Y_4\cup Y_5 =\es$, and hence from \ref{C6-XYZT}, \ref{C6-Y} and \ref{C6-Yi+3}, we conclude that $G[C\cup Y] \in {\cal H}_4$, where we take $S_1:=Y_6$, $S_3:=Y_1$ and $S_2:=Y_2$. This proves \ref{C6-Yi-Yi+1}.
\end{proof2}

   \medskip

   Now if $X\cup Y = \es$, then from \ref{C6-XYZT}, \ref{C6-T} and \ref{C6-T-Z}, we see that $G$ is the join of $G[Z] \cong K_p$ and $G[C\cup T] \cong C_6 + \overline{K_q}$, for some $p\ge 0$ and $q\ge 0$, and we are done; so we may assume that  $X\cup Y \neq  \es$. Then by using \ref{C6-XorY}, we prove the theorem in two cases  as shown below.

   \medskip
\no{\bf Case~1}.~{\it Suppose that $X=\es$ and $Y\neq \es$.}

\smallskip

From \ref{C6-XYZT}, \ref{C6-T} and \ref{C6-T-Z}, we see that $G$ is the join of $H:= G[C\cup Y \cup T] \cong G[C\cup Y] + \overline{K_q}$ and $G[Z] \cong K_p$, for some $p\geq 0$ and $q\ge 0$. So it is enough to show that $ G[C\cup Y]\cong F_2$ or $G[C\cup Y]\in {\cal H}_4 \cup {\cal H}_5$. From \ref{C6-Yi-Yi+1}, we may assume that for each $i\in [6]$, one of $Y_i$ and $Y_{i+1}$ is empty.  Since  $Y\neq \es$, we may assume (up to symmetry) that $Y_6\neq \es$; so we have  $Y_1\cup Y_5=\es$.  Now if $Y_3\neq \es$, then  we have  $Y_2\cup Y_4=\es$, and hence from  \ref{C6-XYZT} and \ref{C6-Y}, we see that $G[C\cup Y] \in {\cal H}_5$ where we take $S_1:=Y_6$ and $S_2:=Y_3$, and we are done; so we may assume that $Y_3=\es$.  Next if $Y_2$ and $Y_4$ are nonempty, then from \ref{C6-Yigeq2}, we conclude that $|Y_6|=1=|Y_2| = |Y_4|$, and hence $G[C\cup Y] \cong F_2$, by \ref{C6-Y}; so we may assume that one of $Y_2$ and $Y_4$ is empty. By symmetry, we may assume that $Y_4=\es$. Then from    \ref{C6-XYZT}, \ref{C6-Y} and \ref{C6-Yi+3}, we see that $G[C\cup Y]\in {\cal H}_4$, where we take $S_1:=Y_6$, $S_2:=Y_2$ and $S_3:=\es$.

  \medskip
  \no{\bf Case 2}.~{\it Suppose that $X\neq \es$ and $Y=\es$.}

\smallskip
 If $X\sm X_i =\es$, for some $i\in \{1,2,3\}$, then  $T = \es$  (by \ref{C6-T-Xi}), and then from \ref{C6-XYZT}, we see that $G$ is the join of $H:=G[C\cup X_i] \in {\cal H}_3$ and $G[Z] \cong K_p$, for some $p\ge 0$, and we are done. So for  each $i\in [3]$, we may assume that  $X\sm X_i \neq \es$. By symmetry and without loss of generality assume that there are vertices, say $x_1\in X_1$ and $x_2\in X_2$. Then by \ref{C6-Xi2}, we have $|X_j|\le 1$ for   $j\in [3]$. Now if there are two vertices in $T$, say $t$ and $t'$, then $\{t,x_1,t',v_2, v_3\}$ induces a $P_2+P_3$ (by \ref{C6-TY} and \ref{C6-TZX}) which is a contradiction; so we may assume that   $|T| \leq 1$. Thus from \ref{C6-Xi} and \ref{C6-TZX}, we conclude that $G[C\cup X\cup T] \in{\cal L}$, and hence  from \ref{C6-XYZT}, \ref{C6-T-Z} and \ref{C6-TZX}, we see that $G$ is the join of $G[C\cup X\cup T]\in {\cal L}$ and $G[Z] \cong K_p$ for some $p\ge 0$.

 This completes the proof. \end{proof}

\section{Some useful results on recoloring} \label{sec:pre-results}

Recall that a $k$-coloring of a graph $G$ is a partition of   $V(G)$ into $k$ independent sets; each such independent set is called a \emph{color class}. Clearly for any graph $G$, we have $\chi(G)\geq \omega(G)$. A \emph{$\chi$-coloring} of a graph $G$ is a $\chi(G)$-coloring of $G$.  Given a coloring $\phi$ of $G$ and a subset $S \subseteq V(G)$, let $\phi(S)$ denote   the set $\{\phi(u)\mid u \in S\}$. For any two colorings $\phi$ and $\psi$ of $G$, we say that $\phi$ = $\psi$ if $\phi(v)$ = $\psi(v)$ for each $v\in V(G)$.  Given two   $k$-colorings of $G$, say $\phi$ and $\psi$,   we say that  $\phi$ and $\psi$ {\em produce the same color classes} if  the set of respective color classes (up to relabelling) produces the same partition of $V(G)$.
Note that a path between any two $k$-colorings of $G$ in   ${\cal R}_{k}(G)$ can be viewed as a {\em recoloring sequence}, where each recoloring step corresponds to the change of color of a single vertex that results in another $k$-coloring of $G$.

We will   use the following known results.
\begin{enumerate}[label=  $(\mathbb{R}\arabic*)$, leftmargin=1cm,series=edu*]
 \item\label{lem:rename}     Bonamy and Bousquet   $\cite{Bonamy2018}$ (see also \cite{Cambie2025}):~({\it  Renaming  Lemma})
    Let $G$ be a graph and let $\ell > \chi(G)$. Let $\phi$ and $\psi$ be  any two $\chi$-colorings of $G$  that produce the same color classes. Then there is a path  between $\phi$ and $\psi$ in ${\cal R}_{\ell}(G)$ where each vertex is recolored at most twice.

\item \label{thm:p3free}    Bonamy and Bousquet  \cite{Bonamy2018}:~Let $G$ be a $P_3$-free graph and let $\ell> \chi(G)$. Then there is a path    between any two $\ell$-colorings of $G$   in ${\cal R}_{\ell}(G)$ where each vertex is recolored at most twice.
\item \label{thm:chordal}   Bonamy and Bousquet  \cite{Bonamy2018}:~Let $G$ be a chordal graph and let $\ell> \chi(G)$. Then  there is a path between any two $\ell$-colorings of $G$  in ${\cal R}_{\ell}(G)$ where each vertex is recolored at most $2n$ times.
\item\label{thm:3k1free}    Merkel  \cite{Owen2022}:~Let $G$ be a $3K_1$-free graph and let $\ell > \chi(G)$. Then there is a path between
    any two $\ell$-colorings of $G$  in ${\cal R}_{\ell}(G)$ where each vertex is recolored at most $4$ times.
 \end{enumerate}

 The following    well known  facts  for the join  and the  union of two graphs, and for a graph with a pair of comparable vertices   have been used by several researchers; see for instance \cite{Manoj-thesis,ref2, Bonamy2018, Cer-Thesis, FJP}. (Below   we assume all the functions on the number of vertices of a given graph are   positive integer valued functions.)
 \begin{enumerate}[label=  $(\mathbb{R}\arabic*)$, leftmargin=1cm,resume=edu*]
 \item\label{lem:join}  Suppose that $G$ is a graph such that for all $\ell_1 > \chi(G)$, there is a path   between any two $\ell_{1}$-colorings of $G$   in ${\cal R}_{\ell_1}(G)$ where each vertex is recolored at most $f(n_G)$ times. Suppose that $H$ is a graph such that for all $\ell_{2} > \chi(H)$, there is a path between any two $\ell_{2}$-colorings of $H$  in ${\cal R}_{\ell_2}(H)$ where each vertex is recolored at most $g(n_H)$ times. Then the following hold: $(i)$~For all $\ell > \chi(G + H)$, there is a path   between any two $\ell$-colorings of $G + H$  in ${\cal R}_{\ell}(G+H)$ where each vertex is recolored at most $\max\{f(n_G), g(n_H)\}$ times. $(ii)$~For all $\ell > \chi(G \vee H)$, there is a path between any two $\ell$-colorings of $G \vee H$ in ${\cal R}_{\ell}(G \vee H)$ where each vertex is recolored at most $\max\{f(n_G), g(n_H)\}$ times.
\item\label{lem:comparable}  Let $G$ be a graph. Suppose that $G$ has a pair of comparable vertices, say $u$ and $v$, such that $N(u)\subseteq N(v)$. Let $H:= G-\{u\}$. Then clearly $\chi(H)=\chi(G)$. Suppose that for all $\ell > \chi(H)$, there is a path between any two $\ell$-colorings of $H$ in ${\cal R}_{\ell}(H)$ where each vertex is recolored at most $f(n_H)$ times, where $f(x)\ge 2x$. Then for all $\ell > \chi(G)$, there is a path between any two $\ell$-colorings of $G$ in ${\cal R}_{\ell}(G)$ where each vertex is recolored at most $f(n_G)$ times.
\end{enumerate}

We will also use the following lemmas.

\begin{lemma}\label{propsition-1}
    Let $G$ be a graph and let $\psi$ be a $\chi$-coloring of $G$. Suppose that for $\ell > \chi(G)$, any $\ell$-coloring $\phi$ of $G$ can be recolored to a $\chi$-coloring $\psi^{*}$ of $G$ by recoloring each vertex of $G$ at most $p$ times where $\psi^{*}$ produces the same color classes as $\psi$. Then $G$ is $\ell$-mixing and the diameter of ${\cal R}_{\ell}(G)$ is at most $2p+2$.
\end{lemma}
\begin{proof}
     Let $\ell > \chi(G)$ and let $\phi_1$ and $\phi_2$ be two $\ell$-colorings of $G$. Then  for each $i\in \{1,2\}$, $\phi_i$ can be recolored to a $\chi$-coloring $\psi_{i}^{*}$   of $G$ by recoloring each vertex of $G$ at most $p$ times where $\psi_i^{*}$   produces the same color classes as $\psi$. Then  by the Renaming Lemma, there is a path  between $\psi_1^{*}$ and $\psi_2^{*}$ in ${\cal R}_{\ell}(G)$ where each vertex of $G$ is recolored at most twice. Thus,  there is a path  between $\phi_1$ and $\phi_2$ in ${\cal R}_{\ell}(G)$ where each vertex of $G$ is recolored at most $2p+2$ times.
\end{proof}

\begin{lemma}\label{prop:independent}
    Let $G$ be a graph and let $\ell > \chi(G)$. Let $I$ be an independent set in $G$ such that $\chi(G)$ = $\chi(G-I)+1$, and let $H :=  G-I$. Suppose that the following hold:
    \begin{itemize}\itemsep=0pt
        \item[(i)]Given any $\ell$-coloring of $G$, one can recolor each vertex of $G$ at most $f(n_G)$ times to obtain an $\ell$-coloring of $G$ where each vertex of $I$ receives the same color that is not on any vertex of $H$.
        \item[(ii)] There is a path between any two $(\ell-1)$-colorings of $H$ in ${\cal R}_{\ell -1}(H)$ where each vertex of $H$ is recolored at most $g(n_{H})$ times.
    \end{itemize}
    Then  there is a path    between any two $\ell$-colorings of $G$   in ${\cal R}_{\ell}(G)$  where each vertex is recolored at most $2[f(n_G) + g(n_{H})]+2$ times.
\end{lemma}
\begin{proof}
    Let $\psi$ be any $\ell$-coloring of $G$. Let $\phi$ be a $\chi$-coloring of $G$ where each vertex in $I$ receives the same color.
    First we prove the following.
  \begin{equation}\label{pr:in}
   \longbox{\em   There is a path between $\psi$ and a $\chi$-coloring $\psi^{*}$ in ${\cal R}_{\ell}(G)$ by recoloring each vertex of $G$ at most $[f(n_G) + g(n_{H})]$ times, where $\psi^{*}$ produces the same color classes as $\phi$.}\tag*{($\star$)}
    \end{equation}
 \no{\em Proof of} \ref{pr:in}:~By $(i)$, starting from $\psi$, we can recolor each vertex of $G$ at most $f(n_G)$ times to obtain an $\ell$-coloring $\beta^{*}$ of $G$ such that $|\beta^*(I)| =1$ and $\beta^*(I) \nsubseteq \beta^*(V(H))$.    Now consider   the restriction of $\beta^{*}$ on $H$, say $\beta^{*}_{H}$. Let $\psi^{*}_{H}$ be an $(\ell-1)$-coloring of $H$ that does not use the color from $\beta^{*}(I)$ and produces the same color classes as $\phi$ when restricted to $H$. From $(ii)$, there is a path between $\beta^{*}_{H}$ and $\psi^{*}_{H}$ in ${\cal R}_{\ell - 1}(H)$, say $\mathcal{P}$,  where each vertex of $H$ is recolored at most $g(n_{H})$ times. Since each coloring in   $\mathcal{P}$ is an $(\ell-1)$-coloring of $H$, we can extend each coloring in  $\mathcal{P}$ to an $\ell$-coloring of $G$ by coloring each vertex in $I$ using the color from $\beta^{*}(I)$. Let $\psi^{*}$ be a $\chi$-coloring of $G$ obtained by extending the coloring $\psi^{*}_{H}$ by  coloring  each vertex in $I$ using the color from  $\beta^{*}(I)$. Therefore, there is a path between $\beta^{*}$ and $\psi^{*}$ in ${\cal R}_{\ell}(G)$   where each vertex of $G$ is recolored at most $g(n_{H})$ times.
    To conclude, there is a path between $\psi$ and $\beta^{*}$ in ${\cal R}_{\ell}(G)$ where each vertex of $G$ is recolored at most $f(n_G)$ times, and there is a path between $\beta^{*}$ and $\psi^{*}$ in ${\cal R}_{\ell}(G)$ where each vertex of $G$ is recolored at most $g(n_H)$ times. Thus, there is a path between $\psi$ and $\psi^{*}$ in ${\cal R}_{\ell}(G)$ where each vertex of $G$ is recolored at most $[f(n_G)+ g(n_H)]$ times.  This proves \ref{pr:in}.  $\sq$

 Now the lemma follows from \ref{pr:in} and \cref{propsition-1}.
        \end{proof}

\section{Recoloring of  ($P_2+P_3, C_4$)-free graphs} \label{sec:recol}

In this section, we prove our main results on recoloring of ($P_2+P_3, C_4$)-free graphs. In fact, our results follow from a sequence of lemmas given below and by using \cref{thm:str-c6-free,case-c6}.


 \begin{lemma}\label{recol-H1}
Let $G\in {\cal H}_1$ and let $\ell> \chi(G)$. Then there is a path between any two $\ell$-colorings of $G$  in ${\cal R}_{\ell}(G)$ where each vertex of $G$ is recolored at most $8$ times.
\end{lemma}
\begin{proof}
Let $G\in {\cal H}_1$  and let $\ell> \chi(G)$. We use the same notation as in the definition of ${\cal H}_1$ (see \cref{gcl2}), and we prove the lemma in two cases.

\smallskip
\no{\bf Case~1}.~$|S|=1$.

Then clearly $\chi(G)=3$  and so $\ell> 3$. We let $S:=\{s\}$.   Let $I:= \{u_3,u_4,u_5\}$ and let $[\ell]$ be the set of colors. Then  $\chi(G-I) =2 = \chi(G)-1$ and clearly $G-I$ is $P_3$-free. Now we will show that given any $\ell$-coloring of $G$, we can recolor each vertex of $G$ at most once to obtain an $\ell$-coloring of $G$ where each vertex of $I$ receives the same color that is not on any vertex of $G-I$. Then the result follows from \ref{thm:p3free} and \cref{prop:independent}.  Let $\phi$ be any $\ell$-coloring of $G$ and  let $\phi(u_1)$ = 1,  $\phi(u_2)$ = 2  and  $\phi(u_5)$ = 3. If $\phi(s)\neq 3$, then we recolor $u_3$  with color 3 and recolor $u_4$ with color 3, and we are done; so we may assume that $\phi(s)=3$. Since $\ell \ge 4$, there is a color, say  $c\in [\ell]\setminus\{1,2,3\}$. If $\phi(u_6) =c$, then starting from $\phi$, we recolor $u_6$ with color 2, and recolor each vertex in $I$, one at a time, with color $c$, and we are done; so $\phi(u_6)\neq c$. Then starting from $\phi$,  we recolor each vertex in $I$, one at a time, with color $c$, and we are done.

\smallskip
\no{\bf Case~2}.~$|S|\geq 2$.

  Let $I:= \{u_3,u_4, u_6\}$ and let $[\ell]$ be the set of colors. Since $G-I$ is $P_3$-free and since $S$ is a clique, $\chi(G-I)  = \max\{|S|, 3\}$. Also since  $S\cup \{u_3\}$ and $\{u_1,u_2,u_5\}$ are  cliques in $G$, we have $\chi(G)\geq \max\{|S|+1, 3\}=|S|+1$  and since     $\chi(G) \leq  \chi(G-I) +1  = \max\{|S|+1, 4\}$. Now if $|S|=2$, then it is easy to see that $\chi(G)=4$, and if $|S|\geq 3$, then $\chi(G)=|S|+1$. In both cases, we have  $\chi(G-I) =\chi(G)-1$.   Now we will show that given any $\ell$-coloring of $G$, we can recolor each vertex of $G$ at most once to obtain an $\ell$-coloring of $G$ where each vertex of $I$ receives the same color that is not on any vertex of $G-I$. Then the result follows from \ref{thm:p3free} and \cref{prop:independent}.
Let $\phi$ be any $\ell$-coloring of $G$ and let $\phi(u_1)$ = 1,  $\phi(u_2)$ = 2  and   $\phi(u_5)$ = 3. Now if there is a color, say $c$, such that  $c\in \phi(I)\sm   \{1, 2, 3\}$, then   $c\notin \phi(S)$; then starting from $\phi$, recolor each vertex in $I$, one at a time, with color $c$ and we are done. So we may assume that $\phi(I)\subseteq \{1, 2, 3\}$. Without loss of generality, assume that $\phi(u_3)$ = 3. Then clearly $3\notin \phi(S)$, and since $\ell \ge 4$, recolor $u_5$ with color 4,  recolor  $u_4$  with color 3, and then recolor $u_6$ with color 3 to get the desired coloring of $G$.
\end{proof}

\begin{lemma}\label{recol-H2}
      Let $G\in {\cal H}_2$ and let $\ell> \chi(G)$. Then there is a path  between any two $\ell$-colorings of $G$ in ${\cal R}_{\ell}(G)$ where each vertex of $G$ is recolored at most $12$ times.
\end{lemma}
\begin{proof}
    Let $G\in {\cal H}_2$ and let $\ell> \chi(G)$. We use the same notation as in the definition of ${\cal H}_2$ (see \cref{gcl2}).  Let $I:= \{u,v,w\}$ and let $[\ell]$ be the set of colors. Then since $V(G-I)$ can be partitioned into two cliques, we see that $G-I$ is   $3K_1$-free; indeed $G-I$ is the complement graph of a bipartite graph and hence perfect; so $\chi(G-I)=\omega(G-I)=\max\{|S_1\cup S_2|,|S_2\cup S_3|, |S_3\cup S_4|\}$. Also since  $S_1 \cup S_2 \cup \{v\}$, $S_2 \cup S_3 \cup \{v\}$ and $S_3 \cup S_4 \cup \{w\}$ are cliques,
     we have $\chi(G)\ge \max\{|S_1\cup S_2|+1,|S_2\cup S_3|+1, |S_3\cup S_4|+1\}  =\chi(G-I)+1$. Now since  $\chi(G)\le \chi(G-I)+\chi(G[I])= \chi(G-I)+1$, we conclude that $\chi(G)=\chi(G-I)+1$.    Then we claim the following.
 \begin{equation} \label{cl-H2}
 \longbox{\em Given any $\ell$-coloring of $G$, we can recolor each vertex of $G$ at most once to obtain an $\ell$-coloring of $G$ where each vertex of $I$ receives the same color that is not on any vertex of $G-I$.}\tag*{($\star$)}
\end{equation}
\no{\em Proof of} \ref{cl-H2}:~
             Let $\phi$ be any $\ell$-coloring of $G$. If $\phi(v)$ = $\phi(w)$, say $\phi(v)$ = $\phi(w)$ = 1, then  $1\notin \phi(S_1\cup S_2 \cup S_3 \cup S_4)$, and starting from $\phi$, we recolor  $u$ with color 1, and we are done. So suppose that $\phi(v) \neq \phi(w)$; let $\phi(v)$ = 1 and   $\phi(w)$ = 2.   If  $1\notin \phi(S_4)$, then starting from $\phi$, we recolor  $u$ with color 1, and recolor $w$ with color 1, and we are done; so we may assume that  $1\in \phi(S_4)$  and let $s_4\in S_4$ be such that $\phi(s_4)=1$. Likewise, we may assume that  $2\in \phi(S_1)$ and let $s_1\in S_1$ be such that $\phi(s_1)=2$.
            Also if there is a color  in $[\ell]$, say $c$, such that  $c\notin \phi(G-\{u\})$, then starting from $\phi$, we recolor $u$, $v$  and $w$, one at a time, with color $c$, and we are done; so we may assume that   $[\ell] \subseteq \phi(G-\{u\})$.
Now since $S_1 \cup S_2 \cup \{v\}$ and $S_3 \cup S_4 \cup \{w\}$ are  cliques in $G$ and since $\ell> \chi(G)\geq \omega(G)$, we have $\ell > |S_1 \cup S_2\cup \{v\}|$ and $\ell > |S_3 \cup S_4 \cup \{w\}|$. Also since $[\ell] \subseteq \phi(G-\{u\})$, there is a color, say $c^*$, such that $c^*\in \phi(S_1 \cup S_2 \cup \{v\})$ or $c^*\in \phi(S_3 \cup S_4 \cup \{w\})$ but not both. If   $c^*\in \phi(S_1 \cup S_2 \cup \{v\})\sm \phi(S_3 \cup S_4 \cup \{w\})$, then starting from $\phi$, we recolor $s_4$ with color $c^*$, recolor  $w$ with color 1, and recolor   $u$ with color 1, and so we have an $\ell$-coloring of $G$ where each vertex in $I$ is colored with 1 and no vertex in $G-I$ is colored with 1, and we are done. Likewise, if   $c^*\in \phi(S_3 \cup S_4 \cup \{w\})\sm \phi(S_1 \cup S_2 \cup \{v\})$, then starting from $\phi$, we recolor $s_1$ with color $c^*$, recolor  $v$ with color 2, and recolor   $u$ with color 2, and so we have an $\ell$-coloring of $G$ where each vertex in $I$ is colored with 2 and no vertex in $G-I$ is colored with 1.  This proves \ref{cl-H2}. $\sq$

Now the lemma follows from \ref{cl-H2}, \ref{thm:3k1free} and \cref{prop:independent}.
\end{proof}

\begin{lemma}\label{thm:recol-C6-free}
    Let $G$ be a $(P_2+P_3, C_4, C_6)$-free graph and let $\ell>\chi(G)$. Then $G$ is $\ell$-mixing and the  diameter of ${\cal R}_{\ell}(G)$ is at most $2n^{2}$.
\end{lemma}
\begin{proof}
    Let $G$ be a $(P_2+P_3, C_4, C_6)$-free graph and let $\ell>\chi(G)$. We will prove the result by showing that there is a path  between any two $\ell$-colorings of $G$ in ${\cal R}_{\ell}(G)$ where each vertex of $G$ is recolored at most $2n$ times. We may assume that $G$ is connected (by \ref{lem:join}), and from \ref{lem:comparable}, we may assume that $G$ does not have a pair of comparable vertices. If $G$ is chordal, then  the result follows from \ref{thm:chordal}. So by \cref{thm:str-c6-free}, we may assume that $G$ is the join of the graph $H$ and a $K_p$, where $H$ is a blowup of a $C_5$ or $H\in \mathcal{H}_1 \cup \mathcal{H}_2$, and $p\ge 0$.  Since a blowup of a $C_5$ is 3$K_1$-free, from  \ref{thm:3k1free} and from \cref{recol-H1,recol-H2}, we conclude that for each $\ell_1 > \chi(H)$,  there is a path between any two $\ell_1$-colorings of $H$ in  ${\cal R}_{{\ell}_1}(H)$  where each vertex of $H$ is recolored at most $2n_H$ times. Since any complete graph is $P_3$-free, the lemma follows from  \ref{thm:p3free} and \ref{lem:join}.
\end{proof}

\begin{lemma}\label{recol-petersen}
    Let $G$ be the Petersen graph and let $\ell > \chi(G)$ = 3. Then there is a path between any two $\ell$-colorings of $G$   in ${\cal R}_{\ell}(G)$ where each vertex is recolored at most 10 times.
\end{lemma}
\begin{proof}
    Let $G$ be the Petersen graph as shown in \cref{gcl2} and let $\ell > \chi(G)$ = 3. Let $I:= \{v_1, v_3, v_9, v_{10}\}$. Then $I$ is a maximum independent set of $G$, and $G-I$ is $P_3$-free.  Now we will show that given any $\ell$-coloring of $G$, we can recolor each vertex of $G$ at most twice to obtain a coloring where each vertex of $I$ receives the same color   that is   not on any vertex of $G-I$. Then the lemma follows from  \ref{thm:p3free} and \cref{prop:independent}.    Let $\phi$ be any $\ell$-coloring of $G$. First we claim the following.
 \begin{equation} \label{cl-pet}
 \longbox{\em If there are at least two vertices in $I$ which receive  the same color under $\phi$,
   then we can recolor each vertex of $G$ at most once to obtain an $\ell$-coloring of $G$ where each vertex of $I$ receives the same color that is not on any vertex of $G-I$.}\tag*{($\star$)}
\end{equation}

\no{\em Proof of} \ref{cl-pet}:~Since the graph is highly symmetric, we may choose any two vertices of $I$, and without loss of generality, we let $\phi(v_9)$ = $\phi(v_{10})$ = 3. If $\phi(v_2) \neq 3$, then starting from $\phi$, we recolor  $v_1$ with color 3 and then recolor $v_3$ with color 3, and we are done; so we may assume that $\phi(v_2)$ = 3. Also if there is a color in $[\ell]$, say $c$, such that $c\notin \phi(N[v_2])$, then starting from $\phi$, we recolor $v_2$ with color $c$, recolor  $v_1$ with color 3, and then recolor $v_3$ with color 3, and we are done. So we may assume that   $[\ell]\subseteq \phi(N[v_2])$, and hence  $\ell =|N[v_2]|= 4$.
    Let $\phi(v_7) = k$, for some $k\in \{1, 2, 4\}$. Then $\phi(v_1) \neq k$ and $\phi(v_3)\neq k$. Now starting from $\phi$, we recolor  $v_7$ with color $\phi(v_1)$, recolor $v_2$ with color $k$,  recolor  $v_1$ with color 3,  and then recolor $v_3$ with color 3, and we are done. Note that we have recolored each vertex of $G$ at most once to reach the desired $\ell$-coloring. This proves \ref{cl-pet}. $\sq$

    \smallskip

By \ref{cl-pet}, we may assume that  no two vertices in $I$ receive  the same color under $\phi$.
    We assume that $\phi(v_1)$ = 1, $\phi(v_3)$ = 2, $\phi(v_9)$ = 3,  and $\phi(v_{10})$ = 4.     Then  $\phi(v_6) \neq 1$, $\phi(v_4) \neq 2$, and $\phi(v_7) \neq 3$. Now starting from $\phi$, we recolor   $v_8$ with color 1, recolor  $v_5$ with color 2, and then recolor   $v_{10}$ with color 3. Now we have a coloring of $G$ where both $v_9$ and $v_{10}$ are colored 3, and  then  we  apply \ref{cl-pet} to get the desired coloring, and we are done.
\end{proof}

\begin{lemma}\label{recol-ind-petersen}
    Let $G \in \cal L$ and let $\ell> \chi(G)$. Then there is a path   between any two $\ell$-colorings of $G$  in ${\cal R}_{\ell}(G)$ where each vertex of $G$ is recolored at most 10 times.
\end{lemma}
\begin{proof}
    Let $G \in \cal L$ and let $\ell> \chi(G)$. Since $G$ contains a $C_5$, we have $\chi(G) = 3$; so $\ell>  3$. Since the Petersen graph is 3-regular, for $\ell\ge 4$, any $\ell$-coloring of $G$ can be extended to an $\ell$-coloring of the Petersen graph. Let $\phi$ and $\psi$ be any two $\ell$-colorings of $G$. Let $\phi^{*}$ and $\psi^{*}$ be two $\ell$-colorings of the Petersen graph obtained by extending $\phi$ and $\psi$, respectively.   Then  by \cref{recol-petersen}, there is a recoloring sequence  between $\phi^{*}$ and $\psi^{*}$, say $\mathcal{P}$,  where each vertex of $G$ is recolored at most 10 times.  Hence, by restricting the recoloring sequence $\mathcal{P}$ to $V(G)$, we obtain a recoloring sequence between $\phi$ and $\psi$ where each vertex of $G$ is recolored at most 10 times.
\end{proof}

\begin{lemma}\label{recol-H3}
    Let $G\in \mathcal{H}_3$ and let $\ell > \chi(G)$. Then there is a path  between any two $\ell$-colorings of $G$  in ${\cal R}_{\ell}(G)$ where each vertex of $G$ is recolored at most 6 times.
\end{lemma}
\begin{proof}
    Let $G\in {\cal H}_3$ and let $\ell> \chi(G)$. We use the same notation as in the definition of ${\cal H}_3$ (see \cref{gcl}).  If $|X|=1$, then $G$ is the theta graph, and the lemma follows from \cref{recol-ind-petersen}. So suppose that  $|X|\geq 2$, and we let $x, x' \in X$. Since $X\cup \{v_1\}$ is a clique, we have $\chi(G)\ge |X|+1\geq 3$.  Let $\phi$ be a coloring of $G$,  where $\phi$ assigns each vertex of $X\cup \{v_1\}$ a distinct color, $\phi(v_4) = \phi(v_1)$, $\phi(v_2) = \phi(x) = \phi(v_6)$  and $\phi(v_3) = \phi(x') = \phi(v_5)$. Then clearly $\phi$ is a ($|X|+1$)-coloring of $G$; so $\chi(G)=|X|+1$  and  $\phi$ is a $\chi$-coloring of $G$. Now we prove that, given any $\ell$-coloring of $G$, we can recolor each vertex of $G$ at most twice to reach a $\chi$-coloring of $G$ that produces the same color classes as  $\phi$. Then the result follows from \cref{propsition-1}.     Let $\psi$ be any $\ell$-coloring of $G$. Suppose $\psi(v_3) \neq \psi(x)$, then we recolor $v_2$ with color $\psi(x)$; otherwise, since $|N[v_3]|< \ell$, there is a color in $[\ell]$, say $c$, such that $c\notin \psi(N[v_3])$, and we recolor $v_3$ with color $c$ and then recolor $v_2$ with color $\psi(x)$. Similarly, we recolor $v_6$ with color $\psi(x)$. Then we recolor  $v_3$ and $v_5$, one at a time, with color $\psi(x')$. Finally, we recolor $v_4$ with color $\psi(v_1)$, and we are done.
\end{proof}

\begin{lemma}\label{recol-H5}
    Let $G\in \mathcal{H}_4\setminus \mathcal{H}_{4}^{*}$ and let $\ell > \chi(G)$. Then there is a path   between any two $\ell$-colorings of $G$  in ${\cal R}_{\ell}(G)$ where each vertex of $G$ is recolored at most 6 times.
\end{lemma}
\begin{proof}
      Let $G\in \mathcal{H}_4\setminus \mathcal{H}_{4}^{*}$  and let $\ell> \chi(G)$. We use the same notation as in the definition of ${\cal H}_4$ (see \cref{gcl}).  We let $S_1:= \{a_1, a_2, \ldots, a_p\}$ and if $S_2$ is nonempty, then we let  $S_2:= \{b_1, b_2, \ldots, b_q\}$. By symmetry, we may assume that $p \ge q$. Since $S_1\cup S_3 \cup \{v_2, v_3\}$ is a clique, we have $\chi(G)\geq |S_1|+|S_3|+2$. Let $\phi$ be a coloring of $G$,  where $\phi$ assigns each vertex of the clique $S_1\cup S_3 \cup \{v_2, v_3\}$ a distinct color,  $\phi(b_i)=\phi(a_i)$, for  $i \in \{1,2,\ldots, q\}$, $\phi(v_1) = \phi(v_5) = \phi(v_3)$  and $\phi(v_4) = \phi(v_6) = \phi(v_2)$. Then clearly $\phi$ is a ($|S_1|+|S_3|+2$)-coloring of $G$; so $\chi(G)=|S_1|+|S_3|+2$ and $\phi$ is a $\chi$-coloring of $G$.    Now we prove that, given any $\ell$-coloring of $G$, we can recolor each vertex of $G$ at most twice to reach a $\chi$-coloring of $G$ that produces the same color classes as  $\phi$. Then the result follows from \cref{propsition-1}. First we show that:
 \begin{equation} \label{cl-H4}
 \mbox{\em  $|S_2| + 3 \leq   \chi(G)< \ell$.}\tag*{($\star$)}
  \end{equation}
\no{\em Proof of} \ref{cl-H4}:~If $S_3\neq\es$, then $|S_2| + 3 \leq  |S_1|+|S_3|+2 = \chi(G)$; so suppose that $S_3=\es$. Then since $G\notin \mathcal{H}_{4}^{*}$, we have $p\neq q$ and hence $|S_1|>|S_2|$. So $|S_2| + 3 \leq |S_1|+2 = \chi(G)$. This proves \ref{cl-H4}. $\sq$

      \medskip
      Now let $\psi$ be any $\ell$-coloring of $G$. If $\psi(v_5)\neq \psi(v_2)$, then we recolor  $v_6$ with   color $\psi(v_2)$; otherwise, since $|N[v_5]| = |S_2| + 3 \leq \chi(G) < \ell$ (by \ref{cl-H4}), there is a color in $[\ell]$, say $c$, such that $c\notin \psi(N[v_5])$, and we  recolor   $v_5$ with color $c$, and then recolor   $v_6$ with color $\psi(v_2)$.  Then we recolor  $v_4$ with color $\psi(v_2)$, recolor  $v_5$ with color $\psi(v_3)$, and then recolor  $v_1$ with color $\psi(v_3)$. Finally, by the Renaming Lemma, we can recolor the vertices in $S_2$ such that color of $b_i$ is the same as the color of $a_i$, for  $i \in \{1,2,\ldots, q\}$ as required, by recoloring each vertex of $S_2$ at most twice, and we are done.
     \end{proof}


\begin{lemma}\label{recol-H62}
    Let $G\in \mathcal{H}_5$  and let $\ell > \chi(G)$. Then there is a path   between any two $\ell$-colorings of $G$ in ${\cal R}_{\ell}(G)$ where each vertex of $G$ is recolored at most 6 times.
\end{lemma}
\begin{proof}
     Let $G\in \mathcal{H}_5$ and let $\ell > \chi(G)$. We use the same notation as in the definition of $\mathcal{H}_5$. We prove the lemma in two cases.

\smallskip
\no{\bf Case~1}.~{\it $|S_1| = 1$ or $|S_2| = 1$.}

By using symmetry, we prove the lemma in this case when $|S_2| = 1$, and we let $S_2:= \{s_2\}$. Since $S_1\neq\es$,   let $s_1 \in S_1$. Since $\{v_6,v_1,v_2,v_3,s_2\}$ induces a $C_5$ and since $S_1$ is a nonempty clique which is complete to $\{v_6,v_1,v_2,v_3,s_2\}$, we have $\chi(G) \geq |S_1| + 3$. Let $\phi$ be a coloring of $G$,  where $\phi$ assigns each vertex of the clique   $S_1\cup  \{v_6,v_1, s_2\}$ a distinct color,   $\phi(v_2) = \phi(s_2)$, $\phi(v_3) = \phi(v_1)$, $\phi(v_4) = \phi(v_6)$ and $\phi(v_5) = \phi(s_1)$. Then clearly $\phi$ is a ($|S_1|+3$)-coloring of $G$; so $\chi(G)=|S_1|+3$ and $\phi$ is a $\chi$-coloring of $G$.  Now we prove that, given any $\ell$-coloring of $G$, we can recolor each vertex of $G$ at most twice to reach a $\chi$-coloring of $G$ that produces the same color classes as  $\phi$. Then the result follows from \cref{propsition-1}. First note that, in any $\ell$-coloring of $G$, for each $j\in \{1,2,4,5\}$, since $|N[v_j]|\le |S_1| + 3  = \chi(G) < \ell$, there is a color in $[\ell]$, say $c_j$,   that does not appear in  $N[v_j]$.     Let $\psi^{*}$ be any $\ell$-coloring of $G$. If $\psi^{*}(v_1) \neq \psi^{*}(s_2)$, then each vertex in $S_1 \cup \{v_6,v_1,s_2\}$ receives a distinct color in $\psi^*$, and if $\psi^{*}(v_1) = \psi^{*}(s_2)$, then we recolor   $v_1$ with color $c_1$. In any case, each vertex in $S_1 \cup \{v_6,v_1,s_2\}$ receives a distinct color in $\psi^*$.   Let this new coloring be $\psi$ (note that $\psi$ may be same as $\psi^{*}$).
Now we recolor $v_2$ with color $\psi(s_2)$. If $\psi(v_4)\neq \psi(v_1)$, then we recolor   $v_3$ with color $\psi(v_1)$; otherwise,  we recolor  $v_4$ with color $c_4$, and then recolor  $v_3$ with color $\psi(v_1)$. Then  recolor  $v_4$ again with color $\psi(v_6)$, and finally recolor   $v_5$ with color $\psi(s_1)$, and we are done.

     \smallskip
\no{\bf Case~2}.~{\it $|S_1| \ge 2$ and $|S_2|\geq 2$.}

We let $s_1, s_1'\in S_1$ and let $s_2, s_2'\in S_2$. Since $S_1 \cup S_2 \cup \{v_6\}$ is a clique in $G$, we have $\chi(G)\geq |S_1|+|S_2|+1$. Let $\phi$ be a coloring of $G$,  where $\phi$ assigns each vertex of the clique $S_1 \cup S_2 \cup \{v_6\}$  a distinct color,    $\phi(v_1) = \phi(s_2)$, $\phi(v_2) = \phi(s_2')$, $\phi(v_3) = \phi(v_6)$, $\phi(v_4) = \phi(s_1')$  and $\phi(v_5) = \phi(s_1)$. Then clearly $\phi$ is a ($|S_1|+|S_2|+1$)-coloring of $G$; so $\chi(G)=|S_1|+|S_2|+1$ and $\phi$ is a $\chi$-coloring of $G$.    Now we prove that, given any $\ell$-coloring of $G$, we can recolor each vertex of $G$ at most twice to reach a $\chi$-coloring of $G$ that produces the same color classes as  $\phi$. Then the result follows from \cref{propsition-1}. First note that, for each $j\in \{1,2,4,5\}$, since $|N[v_j]|\le |S_1|+|S_2|+1  = \chi(G) < \ell$, in any $\ell$-coloring of $G$, there is a color in $[\ell]$, say $c_j$, that does not appear in $N[v_j]$.
    Let $\psi$ be any $\ell$-coloring of $G$. If $\psi(v_2) \neq \psi(s_2)$, then recolor  $v_1$ with color $\psi(s_2)$; otherwise, we recolor  $v_2$ with color $c_2$, and then recolor $v_1$ with color $\psi(s_2)$. Likewise, we recolor   $v_5$ with color $\psi(s_1)$. Since  $\{v_3\}$ is complete to $S_1\cup S_2$, we have $\psi(s_2') \neq \psi(v_3) \neq \psi(s_1')$. Now we recolor  $v_2$ with color $\psi(s_2')$, recolor  $v_4$ with color $\psi(s_1')$, and recolor $v_3$ with color $\psi(v_6)$, and we are done.
\end{proof}

  \begin{lemma}\label{H-frozen}
If $G\in \{C_6, F_2\}\cup  {\cal H}_4^{*}$, then  $G$ admits a frozen ($\chi(G)+1$)-coloring.
\end{lemma}
\begin{proof}
First note that $\chi(C_6)=2$  and  $\chi(F_2)=3$. Now since the graph $C_6$ admits a frozen $3$-coloring  (see \cref{cat}:$(ii)$) and since  the graph $F_2$ admits a frozen $4$-coloring  (see \cref{cat}:$(iv)$), it is enough to prove the lemma when $G\in {\cal H}_4^*$. We use the same notation as in the definition of ${\cal H}_4^*$.
 Since $S_1\cup \{v_1,v_2\}$ is a clique in $G$, we have $\chi(G)\geq p+2$. Define $\phi: V(G)\rightarrow \{1,2,\ldots, p+2\}$ where
 $\phi(v_1)=1 =\phi(v_3)=\phi(v_5)$,  $\phi(v_2)=2 =\phi(v_4)=\phi(v_6)$ and $\phi(S_1)=\phi(S_2)=\{3,4,\ldots, p+2\}$.
Then clearly $\phi$ is a ($p+2$)-coloring of $G$, and hence $\chi(G)= p+2$. If $p=1$, then $G\cong F_1$, and see  \cref{cat}:$(iii)$ for a frozen $4$-coloring of $G$; so assume that $p\geq 2$. Now we define $\psi: V(G)\rightarrow \{1,2,\ldots, p+3\}$  where
$\psi(v_1)=\psi(v_4) =1$, $\psi(v_2)=2$, $\psi(v_5)=2$, $\psi(v_3)= \psi(v_6) =3$, $\psi(S_1)=\{4,5,\ldots, p+3\}$ and  $\psi(S_2)=\{4,5,\ldots, p+3\}$.
Then it is easy to verify that $\psi$ is a frozen ($p+3$)-coloring of $G$, and we are done.
\end{proof}

\begin{lemma}\label{union-join}
If $H\in \{C_6, F_2\}\cup  {\cal H}_4^{*}$, then for $p\geq 0$ and $q\geq 0$,   $(H+\overline{K_q})\vee K_p$ is not recolorable.
\end{lemma}
\begin{proof}
It is known that (see \cite{Manoj-thesis}) for any two graphs $G_1$ and $G_2$, the graph $G_1+G_2$ (resp. $G_1\vee G_2$) is recolorable if and only if $G_1$ and $G_2$ are recolorable.
  Now if $H\in \{C_6, F_2\}\cup  {\cal H}_4^{*}$, then from \cref{H-frozen}, we see that $H$ is not recolorable.  Hence $(H+\overline{K_q})\vee K_p$ is not recolorable.
\end{proof}

 Now we give a complete characterization of   $(P_2+P_3, C_4)$-free graphs that are recolorable.

\begin{theorem}\label{recol-main}
Let $G$ be a connected ($P_2+P_3, C_4$)-free graph. Then $G$ is recolorable   if and only if $G$ is not isomorphic to $(H+\overline{K_q})\vee K_p$, where  $H\in \{C_6, F_2\}\cup  {\cal H}_4^{*}$, $p\ge 0$ and $q\geq 0$. Moreover, if $G$ is recolorable, then for all $\ell > \chi(G)$,  the  diameter of ${\cal R}_{\ell}(G)$ is at most 2$n^{2}$.
\end{theorem}
\begin{proof} Let $G$ be a connected ($P_2+P_3, C_4$)-free graph. If $G$ is  isomorphic to $(H+\overline{K_q})\vee K_p$, where  $H\in \{C_6, F_2\}\cup  {\cal H}_4^{*}$, $p\ge 0$ and $q\geq 0$, then $G$ is not recolorable, by  \cref{union-join}. Conversely, suppose that $G$ is not isomorphic to $(H+\overline{K_q})\vee K_p$, where  $H\in \{C_6, F_2\}\cup  {\cal H}_4^{*}$, $p\ge 0$ and $q\geq 0$.  Now it is enough show that for all $\ell> \chi(G)$, $G$ is $\ell$-mixing and the diameter of ${\cal R}_{\ell}(G)$ is at most 2$n^{2}$.  We will prove this by showing that there is a path  between any two $\ell$-colorings of $G$ in ${\cal R}_{\ell}(G)$ where each vertex of $G$ is recolored at most $2n$ times.
 By \cref{thm:recol-C6-free}, we may assume that $G$ contains a $C_6$, and from \ref{lem:comparable}, we may assume that $G$ does not have a pair of  comparable vertices. Then by \cref{case-c6}, $G$ is the join of the graph $(H + \overline{K_q})$ and a $K_p$, where $H \in {\cal L} \cup \mathcal{H}_3 \cup (\mathcal{H}_4\setminus \mathcal{H}_4^{*}) \cup \mathcal{H}_5$, $p\geq 0$ and $q\ge 0$. Note that, by \ref{thm:p3free}, for all $\ell_1 > p$, there is a path between any two $\ell_1$-colorings of $K_p$ in ${\cal R}_{{\ell}_1}(K_p)$ where each vertex of $K_p$ is recolored at most twice.
Then the theorem follows from \cref{recol-ind-petersen,recol-H3,recol-H5,recol-H62} and  from \ref{lem:join}. \end{proof}

Next we will  prove   that   \cref{Cer-conj} holds for the class of ($P_2+P_3, C_4$)-free graphs. We   need the following results.

\begin{lemma}[\cite{ref2}]\label{recol-cycle}
    For all $t\geq3$ and $k \ge 4$,    $C_t$  is $k$-mixing and the  diameter of ${\cal R}_{k}(C_t)$  is at most $4n$.
\end{lemma}

\begin{lemma}\label{recol-F2}
   For all $k\ge \rho(F_2)+2$,  $F_2$ is $k$-mixing and the  diameter of ${\cal R}_{k}(F_2)$  is at most $6n$.
\end{lemma}
\begin{proof}
    Let $F_2$ be the graph as shown in   \cref{cat}. Since $\{s_1, v_1, v_2\}$ is a clique, we have $\chi(G)\geq 3$.  Let $\phi$ be a coloring of $G$,  where $\phi$ assigns each vertex of the clique $\{s_1, v_1, v_2\}$ a distinct color, $\phi(v_3) = \phi(v_5) = \phi(v_1)$, $\phi(v_4) = \phi(v_6) = \phi(v_2)$, and $\phi(s_2) = \phi(s_3) = \phi(s_1)$. Then clearly $\phi$ is a $3$-coloring of $G$; so $\chi(G)=3$ and $\phi$ is a $\chi$-coloring of $G$.  Now we prove that, given any $k$-coloring of $G$, we can recolor each vertex of $G$ at most twice to reach a $\chi$-coloring of $G$ that produces the same color classes as  $\phi$. Then the result follows from  \cref{propsition-1}.
     Let $\psi$ be an $k$-coloring of $F_2$. If $\psi(v_4)\neq \psi(v_1) \neq \psi(s_2)$, then we recolor $v_3$ with color $\psi(v_1)$; otherwise,  if $\psi(v_4) = \psi(v_1)$ (say), then since $|N[v_4]| = 5 < k$, there is a color in $[k]$, say $c$, such that $c\notin \psi(N[v_4])$, and we recolor $v_4$ with color $c$, and then recolor $v_3$ with color $\psi(v_1)$. Similarly, if $\psi(s_2) = \psi(v_1)$, then since $|N[s_2]| = 5 < k$, there is a color in $[k]$, say $c^*$, such that $c^* \notin \psi(N[s_2])$, and we recolor $s_2$ with color $c^*$ and then recolor $v_3$ with color $\psi(v_1)$. Likewise, similar to $v_3$, we can recolor $v_6$ with color $\psi(v_2)$. Then recolor $v_4$ with color $\psi(v_2)$, and recolor $v_5$ with color $\psi(v_1)$. Finally, we recolor  $s_2$ and $s_3$, one at a time, with color $\psi(s_1)$, and we are done.
\end{proof}

We observe that the proof of \cref{recol-H5} is based on the fact that $|N[v_5]| \le \chi(G) < k$ and hence, in  any $k$-coloring of $G$,  we always have an extra color  that does not appear in $N[v_5]$. For a graph $G \in \mathcal{H}_{4}^{*}$, if $k \ge \rho(G) +2$, then since $\rho(G) \ge |N(v_5)|$, we have $|N[v_5]| < k$ and hence, in  any $k$-coloring of $G$,  we always have a  color  that does not appear in $N[v_5]$. So we have the following.

\begin{lemma}\label{recol-H4*}
 If $G\in \mathcal{H}_4^{*}$, then for $k\ge \rho(G)+2$, $G$ is $k$-mixing and the  diameter of ${\cal R}_{k}(G)$  is at most $6n$.
\end{lemma}

\begin{theorem}\label{thm:degen}
    Let $G$ be a $(P_2+P_3, C_4)$-free graph. Then, for all $k\ge \rho(G)+2$,   the   diameter of ${\cal R}_{k}(G)$ is at most $2n^{2}$.
\end{theorem}
\begin{proof}Let $G$ be a $(P_2+P_3, C_4)$-free graph and let $k\ge \rho(G)+2$. By \ref{lem:join}:$(ii)$, we may assume that $G$ is connected. Recall that for any graph $G^*$, we have  $\chi(G^*) \le \rho(G^*)+1$, and hence $\chi(G^*)+1 \le \rho(G^*)+2$. Therefore, by \cref{recol-main}, we may assume that $G$ is  isomorphic to $(H+\overline{K_q})\vee K_p$, where  $H\in \{C_6, F_2\}\cup  {\cal H}_4^{*}$, $p\ge 0$ and $q\geq 0$. Note that, by \ref{thm:p3free}, for all $k_1 > p$, there is a path between any two $k_1$-colorings of $K_p$ in ${\cal R}_{k_1}(K_p)$ where each vertex of $K_p$ is recolored at most twice. Now the theorem follows from \cref{recol-cycle,recol-F2,recol-H4*}, and from \ref{lem:join}. \end{proof}

\bigskip
\no{\bf Acknowledgement}.~The first author  would like to acknowledge  the financial support from Anusandhan National Research Foundation (ANRF), Government of India,  under the National Post Doctoral Fellowship (N-PDF) program (No.~PDF/2025/005208), and by the Natural Sciences and Engineering Research Council of Canada (NSERC) grant RGPIN-2016-06517.

\end{document}